\documentclass[]{interact}

\usepackage{epstopdf}
\usepackage[caption=false]{subfig}

\theoremstyle{plain}
\newtheorem{theo}{Theorem}[section]
\newtheorem{lem}[theo]{Lemma}
\newtheorem{cor}[theo]{Corollary}
\newtheorem{prop}[theo]{Proposition}

\theoremstyle{definition}
\newtheorem{defin}[theo]{Definition}
\newtheorem{example}[theo]{Example}

\theoremstyle{remark}
\newtheorem{rem}{Remark}
\newtheorem*{notation*}{Notation}

\newcommand\C{\mathbb{C}}
\newcommand\Z{\mathbb{Z}}

\newcommand\R{\mathbb{R}}

\newcommand\g{\mathfrak{g}}

\newcommand\p{\mathfrak{p}}
\newcommand\h{\mathfrak{h}}
\newcommand\n{\mathfrak{n}}

\newcommand\fb{\mathfrak{b}}
\newcommand\fsl{\mathfrak{sl}}
\newcommand\fsp{\mathfrak{sp}}
\newcommand\fo{\mathfrak{o}}
\newcommand\fgl{\mathfrak{gl}}
\newcommand\fl{\mathfrak{l}}

\newcommand\fS{\mathfrak{S}}

\newcommand\fu{\mathfrak{u}}
\newcommand\fk{\mathfrak{k}}

\newcommand\sm{\mathsf{m}}

\newcommand\D{\Delta}

\newcommand\cF{\mathcal{F}}

\newcommand\cI{\mathcal{I}}

\newcommand\cS{\mathcal{S}}

\newcommand\cD{\mathcal{D}}
\newcommand\cO{\mathcal{O}}

\newcommand\bx{\mathbf{x}}

\newcommand\bU{\mathbf{U}}

\DeclareMathOperator{\Int}{Int}

\DeclareMathOperator{\Span}{Span}

\DeclareMathOperator{\Supp}{Supp} 

\begin{document}


\title{On the structure of simple bounded weight modules of $\fsl(\infty)$, $\fo(\infty)$, $\fsp(\infty)$}

\author{\name{L. Calixto\thanks{Email: lhcalixto@ufmg.br}}
\affil{Department of Mathematics, Federal University of Minas Gerais, Belo Horizonte, Brazil}
}

\maketitle

\begin{abstract}
We study the structure of bounded simple weight $\fsl(\infty)$-, $\fo(\infty)$-, $\fsp(\infty)$-modules, which have been recently classified in \cite{GP18}. Given a splitting parabolic subalgebra $\p$ of $\fsl(\infty)$, $\fo(\infty)$, $\fsp(\infty)$, we introduce the concepts of $\p$-aligned and pseudo $\p$-aligned \mbox{$\fsl(\infty)$-,}$\fo(\infty)$-, $\fsp(\infty)$-modules, and give necessary and sufficient conditions for bounded simple weight modules to be $\p$-aligned or pseudo $\p$-aligned. The existence of pseudo $\p$-aligned modules is a consequence of the fact that the Lie algebras considered have infinite rank.
\end{abstract}


\begin{amscode}
17B65, 17B10
\end{amscode}

\section{Introduction}

In the representation theory of finite-dimensinal Lie algebras it is very important to answer the following natural questions for a simple $\g$-module $M$:
\begin{enumerate}
\item \label{Q1} given a Borel subalgebra $\fb$ of $\g$, is $M$ a $\fb$-highest weight module?
\item \label{Q2} given a parabolic subalgebra $\p$ of $\g$, is $M$ an irreducible quotient of a generalized Verma module associated to $\p$?
\end{enumerate}
These concepts lead to the study of well known categories, namely, the category $\cO$ and the parabolic categories $\cO$. Modules which do not belong to these categories can also be of great interest, for instance, Harish-Chandra modules (see \cite{Dix96}).

Weight modules for Lie algebras have been extensively studied since 1970's by several authors: G. Benkart, D. Britten, S. Fernando, V. Futorny, A. Joseph, F. Lemire,  V. Mazorchuk, I. Penkov, V. Serganova and others. The classification of simple modules with finite-dimensional weight spaces over a reductive finite-dimensional Lie algebra was completed by O. Mathieu in his remarkable work \cite{Mat00}. This classification relies on the parabolic induction theorem of Fernando and Futorny, \cite{Fer90, Fut87}, which states that any weight module with finite dimensional weight spaces is isomorphic to the irreducible quotient of a certain parabolically induced module. 

To answer the two questions above in the finite-dimensional case, in general it is enough to understand the Fernando-Kac subalgebra associated to the module $M$; for the definition of this subalgebra see for instance \cite[Page~2]{PSZ04}, but also \cite{Fer90, Kac85}. The Fernando-Kac subalgebra is in turn closely related to the notion of the shadow of a module, which was introduced by I. Dimitrov, O. Mathieu and I. Penkov in \cite{DMP04}. In \cite{DP99}, this notion was extended to the case of infinite-dimensional root reductive Lie algebras. In particular, it was shown that unlike the finite-dimensional case, the shadow does not give sufficient information to anwer questions \eqref{Q1} and \eqref{Q2}. For instance, it may happen that the shadow of $M$ coincides with the shadow of a highest weight module without $M$ being a highest weight module. In the recent paper \cite{GP18}, D. Grantcharov and I. Penkov  classify simple bounded weight modules over the infinite-dimensional Lie algebras $\fsl(\infty)$, $\fo(\infty)$ and $\fsp(\infty)$. There the authors also discuss when a given module is a highest weight module.

This brings us to the natural problem, which is the main goal of the current paper, of answering question \eqref{Q2} for the modules that appear in the classification given in \cite{GP18}. Let $\g$ be equal $\fsl(\infty)$, $\fsp(\infty)$. It is standard to represent $\g$ as a direct limit Lie algebra $\g=\varinjlim_n \g_n$, where $\g_n$ are simple classical Lie algebras of the same type and ${\rm rank}\ \g_n=n$; if $\g=\fo(\infty)$, then all $\g_n$ are either of type $B$ or $D$. Let $\p$ be a splitting parabolic subalgebra of $\g$, that is, a subalgebra such that the intersections $\p_n:=\p\cap \g_n$ are parabolic subalgebras of $\g_n$. In the current paper, we call a simple weight module $M$ \emph{$\p$-aligned} if it is isomorphic to the irreducible quotient of a generalized Verma module associated to $\p$. If $M$ is not $\p$-aligned but $M\cong \varinjlim_n M_n$, where $M_n$ is a simple $\p_n$-aligned module for every $n$ (i.e. $M_n$ is a quotient of a generalized Verma module associated to $\p_n$),  we call the module $M$ \emph{pseudo $\p$-aligned}. 

Our main result is an explicit criterion for each simple bounded weight module $M$ from the classification in \cite{GP18}, and for each splitting parabolic subalgebra $\p$ containing a fixed Cartan subalgebra, to be $\p$-aligned or pseudo $\p$-aligned. We also compute the relevant inducing $\p$-module in the case when $M$ is $\p$-aligned. The results are presented in two sections, where we consider separately the cases that $M$ is integrable and non-integrable.

\begin{notation*}
In this paper we follow the notation of \cite{GP18}. Namely, the ground field is $\C$, and all vector spaces, algebras, and tensor products are considered to be over $\C$, unless otherwise declared. We let $\langle \ \rangle_\C$ denote the span over $\C$. We write $\C^{\Z_{>0}}$ for the set of all sequences of complex numbers, and $\C_{\rm f}^{\Z_{>0}}$ for the subset of $\C^{\Z_{>0}}$ consisting all sequences that admit only finitely many nonzero entries. The symmetric and exterior powers of a vector space will be denoted by $S^\cdot(\cdot)$ and $\Lambda^\cdot(\cdot)$, respectively. The universal enveloping algebra of a Lie algebra $\fk$ will be denoted by $\bU(\fk)$.
\end{notation*}

\section{\bf{Preliminaries}}

\subsection{Generalities}\label{subsec:generalities} Throughout the paper $\g_n$ will be $\fsl(n+1)$, $\fo(2n+1)$, $\fo(2n)$, $\fsp(2n)$. A fixed Cartan subalgebra of $\g_n$ will be denoted by $\h_n$, and we let $\varepsilon_i$ denote the standard vectors of $\h_n^*$ (using Bourbaki's notation). For $\fo(2n+1)$, $\fo(2n)$, $\fsp(2n)$ we identify $\h_n^*$ with $\C^n$: any $(\lambda_1,\ldots, \lambda_n)\in \C^n$ is identified with the element $\sum \lambda_i\varepsilon_i\in \h_n^*$. Similarly, we  identify $\C^{n+1}$ with the Cartan subalgebra of $\fgl(n+1)$ and when we say that $\lambda\in \C^{n+1}$ is a weight of $\fsl(n+1)$, we are considering the projection of $\lambda$ in $\h_n^*\cong \C^{n+1}/\C {\bf 1}$, where ${\bf 1} = (1,\ldots, 1)\in \C^{n+1}$.

Let $M$ be a $\g_n$-module. We call $M$ an \emph{integrable module} if for any $m\in M$, $g\in \g_n$, we have $\dim \langle m,\ g\cdot m,\ g^2\cdot m,\ldots\rangle_\C<\infty$. We call $M$ a \emph{weight module} if it admits a weight space decomposition: $M=\bigoplus_{\mu\in \h_n^*} M_\mu$, where $M_\mu=\{m\in M\mid hm=\mu(h)m,\ \forall h\in \h_n^*\}$. The \emph{support} of a weight module $M$ is the set $\Supp M := \{\mu\in \h_n^*\mid M_\mu\neq 0\}$. All modules considered in this paper are assumed to be weight modules with finite-dimensional weight spaces, that is, $\dim M_\mu < \infty$ for every $\mu\in \Supp M$. A weight module $M$ is called \emph{bounded} if there is a constant $c\in \Z_{>0}$ for which $\dim M_\mu < c$ for all $\mu\in \Supp M$. If $\fk$ is a subalgebra of $\g_n$, then a weight $\nu\in \Supp M$ is called \emph{$\fk$-singular} if and only if $\fk\cdot M_\nu=0$. A nonzero vector $v\in M$ is $\fk$-singular if $\fk\cdot v=0$. The Lie algebras $\fsl(2n+1)$, $\fo(2n+1)$, $\fo(2n)$, $\fsp(2n)$ admit a \emph{natural representation} that will be denoted by $V_{n+1}$, $V_{2n+1}$, $V_{2n}$, $V_{2n}$, respectively. These representations are characterized, up to isomorphism, by their supports: $\Supp V_{n+1}=\{\varepsilon_i\mid 1\leq i\leq n+1\}$ for $\g_n = \fsl(n+1)$; $\Supp V_{2n+1}=\{0, \pm \varepsilon_i\mid 1\leq i\leq n\}$ for $\g_n = \fo(2n+1)$; and $\Supp V_{2n}=\{\pm \varepsilon_i\mid 1\leq i\leq n\}$ for $\g_n = \fo(2n),\ \fsp(2n)$.

Throughout the paper $\g$ will be equal $\fsl(\infty)$, $\fsp(\infty)$, $\fo(\infty)$. Then $\g=\varinjlim_n \g_n$, where each $\g_n$ is a simple classical Lie algebra of the same type with ${\rm rank}\ \g_n=n$. The embeddings defining $\g$ are always assumed to be root embeddings, that is, $\h_n$ is mapped into $\h_{n+1}$, and any root space of $\g_n$ is mapped into a root space of $\g_{n+1}$. From now on we fix $\h:=\varinjlim_n \h_n$. Then $\h$ is a \emph{Cartan subalgebra} of $\g$, and the adjoint action of $\h$ on $\g$ yields a root space decomposition $\g = \h\oplus \left(\bigoplus_{\alpha\in \h^*}\g_\alpha\right)$, where $\g_\alpha = \{x\in \g\mid [h,x]=\alpha(h)x,\ \forall h\in \h\}$. The set of roots of $\g$ is $\D = \{\alpha\in \h^*\setminus \{0\}\mid \g_\alpha\neq 0\}$.  When needed, we write $\h_A$, $\h_B$, $\h_C$, $\h_D$ to let clear which family of embeddings we are considering, and we let $\D_A$, $\D_B$, $\D_C$, $\D_D$ denote the root systems relative to each of these Cartan subalgebras, respectively. Explicitly we have that
\begin{gather*}
\D_A=\{\varepsilon_i-\varepsilon_j\mid i\neq j\};\quad \D_C=\{\pm 2\varepsilon_i,\pm\varepsilon_i \pm\varepsilon_j\mid i\neq j\};  \\
\D_B=\{\pm\varepsilon_i,\pm\varepsilon_i \pm\varepsilon_j\mid i\neq j\};\quad \D_D=\{\pm\varepsilon_i \pm\varepsilon_j\mid i\neq j\},
\end{gather*}
where $i,j\in \Z_{>0}$, and $\varepsilon_i$ stands by the vectors of $\h^*$ whose restriction to $\h_n^*$ coincide with the vectors $\varepsilon_i\in \h_n^*$ defined above for every $n$. We identify an element $\lambda\in \h^*$ with the formal sum $\sum \lambda_i\varepsilon_i$, or with the infinite sequence $(\lambda_1,\lambda_2, \ldots)\in \C^{\Z_{>0}}$. Notice that if $\g=\fsl(\infty)$, then two infinite sequences in $\C^{\Z_{>0}}$ determine the same element in $\h^*$ if their difference is a constant sequence ${\bf c}:=(c,c,\ldots)\in \C^{\Z_{>0}}$, for some $c\in \C$. In what follows we let $Q_\g:=\sum_{\alpha\in \D} \Z\alpha$ be the root lattice of $\g$ associated to $\h$.

Every Borel subalgebra $\fb$ of $\g$ we consider in this paper is assumed to be a \emph{splitting Borel subalgebra} containing the fixed Cartan subalgebra $\h$, that is, all the intersections $\fb_n:=\fb\cap \g_n$ are Borel subalgebras of $\g_n$ which contain $\h_n$. In particular, $\fb=\h\oplus \left(\bigoplus_{\alpha\in \D^+} \g_\alpha\right)$, for some triangular decomposition $\D=\D^-\cup \D^+$. The terminology introduced above for $\g_n$-modules also make sense for $\g$-modules and it will be used freely. Moreover, we let $L_\fb(\lambda)$ denote a $\fb$-highest weight $\g$-module with highest weight $\lambda\in \h^*$. Similarly we define a $\g_n$-module $L_{\fb_n}(\lambda^n)$ for any $\lambda^n\in \h_n^*$.

\subsection{Parabolic subalgebras}
In this paper, a subalgebra $\p$ of $\g$ is called \emph{parabolic} if it contains some Borel subalgebra of $\g$. In particular, according to our definition of Borel subalgebras, we must have $\h\subseteq \p$. Parabolic subalgebras are in bijection with parabolic sets $P\subseteq \D$ (i.e. $P$ is additively closed and $\D=\{-P\}\cup P$). Namely, given a parabolic subalgebra $\p$, we define $P_\p=\{\alpha\in \D\mid \g_\alpha\subseteq \p\}$, and, on the other hand, for  a given parabolic set $P$, we define $\p(P)=\h\oplus \bigoplus_{\alpha\in P}\g_\alpha$. If we set $P^0:=P\cap\{-P\}$ and $P^+=P\setminus P^0$, then  $\fl(P):=\h\oplus \bigoplus_{\alpha\in P^0}\g_\alpha$ is the \emph{reductive component} of $\p(P)$ while $\fu(P):=\bigoplus_{\alpha\in P^+}\g_\alpha$ is its locally \emph{nilpotent radical}.

Let $(F,\prec)$ be a partially ordered set, and, for every $i\in F$, define 
	\[
[i]_\prec :=\{i\}\cup \{j\in F\mid i\text{ and }j\text{ are not comparable}\}.
	\]
When there is no risk of confusion, we write $[i]$ instead of $[i]_\prec$. If we say $i\approx j$ if and only if $[i]=[j]$, then $\approx$ defines an equivalence relation on $F$. Let $[F]$ denote the set of equivalence classes $F/\approx$.

\begin{defin}
A partial order $\prec$ on a set $F$ is called \emph{admissible} if the following condition is satisfied: $[i]\neq [j]$ implies $[i]\cap [j]=\emptyset$, for all $i,j\in F$ (this is equivalent to saying that $j\in [i]$ implies $[j]=[i]$, for all $i,j\in F$).
\end{defin}

It is easy to see that $\prec$ induces a linear order on $[F]$ if and only if $\prec$ is admissible. In the current paper we will always assume that a partial order is admissible. If $A, B\subseteq F$, by writing $A\prec B$ we mean that $a\prec b$ for any $a\in A$ and $b\in B$.

A partial order $\prec$ on $\Z$ is called \emph{$\Z_2$-linear} if multiplication by $-1$ reverses the order. Let $\prec$ be a partial order on $\Z_{>0}$ (resp. a $\Z_2$-linear partial order on $\Z$). Then $\prec$ induces a partial order on the set $\{\varepsilon_i\}$ (resp. a $\Z_2$-linear partial order on the set $\{\pm \varepsilon_i\}\cup \{0\}$), which in turn induces a partial (resp. a $\Z_2$-linear partial) order on $\D\cup \{0\}$. Explicitly we have 
\begin{equation}\label{eq:order.Delta=order.Z}
\begin{aligned}
\varepsilon_i - \varepsilon_j \succ 0 \Leftrightarrow \varepsilon_i \succ \varepsilon_j \Leftrightarrow i\prec j, \\
\varepsilon_i + \varepsilon_j \succ 0 \Leftrightarrow \varepsilon_i \succ -\varepsilon_j \Leftrightarrow i\prec -j.
\end{aligned}
\end{equation}
Define $P(\prec) = \{\alpha\in \D\mid \alpha\succ 0\text{ or }\alpha \text{ not comparable with }0\}$, and let $\p(\prec)$ be the subalgebra of $\g$ generated by $\h$ and by all root spaces $\g_\alpha$ with $\alpha\in P(\prec)$. In the other direction, for each parabolic subalgebra $\p\subseteq \g$ we define a partial (resp. a $\Z_2$-linear partial) order $\prec$ on $\D\cup \{0\}$ by setting $\alpha\succ 0$ if and only if $\g_\alpha\subseteq \p$ but $\g_{-\alpha}\nsubseteq \p$. Using \eqref{eq:order.Delta=order.Z} it is clear that $\prec$ also determines a partial order on $\Z_{>0}$ (resp. a $\Z_2$-linear partial order on $\Z$).

\begin{lem}\cite[Proposition~4]{DP99}\label{lem:p.ord.<-->par.subalgebras}
The following statements hold:
\begin{enumerate}
\item \label{lem:p.ord.<-->par.subalgebras.A} If $\g=\fsl(\infty)$, then there is a bijection between parabolic subalgebras of $\g$ and admissible partial orders on $\Z_{>0}$.
\item  \label{lem:p.ord.<-->par.subalgebras.B,C} If $\g=B(\infty)$, or $C(\infty)$, then there is a bijection between  parabolic subalgebras of $\g$ and admissible $\Z_2$-linear partial orders on $\Z$.
\item If $\g=D(\infty)$, then there is a bijection between parabolic subalgebras of $\g$ and admissible $\Z_2$-linear partial orders on $\Z$ with the property that if $i\in \Z\setminus\{0\}$ is not comparable with $0$, then $j$ is also not comparable with $0$ for some $j\neq i$.
\end{enumerate}
\end{lem}

From now on we always assume that partial orders are admissible. Let $\prec$ be a ($\Z_2$-linear) partial order on $\Z_{>0}$ (resp. $\Z$), and let $\p(\prec)$ denote the parabolic subalgebra of $\g$ associated to $\prec$. Let $\fl(\prec)$ denote the reductive component of $\p(\prec)$, and $\fu(\prec)$ denote the locally nilpotent radical of $\p(\prec)$. Set $\fl(\prec)^{ss} := [\fl(\prec), \fl(\prec)]$. Notice that $\fl(\prec)^{ss}$ is generated by all root spaces $\g_\alpha$ with $\alpha\in P(\prec)^0 = \{\alpha\in \D\mid \alpha \text{ not comparable with }0\}$, and $\fu(\prec)$ is generated by all root spaces $\g_\alpha$ with $\alpha\in P(\prec)^+ = \{\alpha\in \D\mid \alpha\succ 0\}$.

\begin{rem}\label{rem:reductive.nilp.comp}
\begin{enumerate}
\item Let $\g=A(\infty)$. Then we have a decomposition $P(\prec)^0=\cup P(\prec)_{[p]}^0$, where $P(\prec)_{[p]}^0=\{\varepsilon_i-\varepsilon_j\mid i,j\in [p],\ i\neq j\}$. If $\g[p]$ denotes the Lie subalgebra of $\g$ generated by the root spaces $\g_{\alpha}$ for $\alpha \in P(\prec)_{[p]}^0$, then $\g[p]$ is of type $ A$. Moreover, $P(\prec)^+=\{\varepsilon_i-\varepsilon_j\mid i\prec j\}$.
\item If $\g=D(\infty)$, then notice that: $\pm(\varepsilon_i-\varepsilon_j)\in P(\prec)^0$ if and only if $[i]=[j]$; $\pm(\varepsilon_i + \varepsilon_j)\in P(\prec)^0$ if and only if $[i]=[-j]$; and $\{\pm(\varepsilon_i+\varepsilon_j), \ \pm(\varepsilon_i-\varepsilon_j)\}\subseteq P(\prec)^0$ if and only if $[i]=[j]= [0]$. Thus we have a decomposition $P(\prec)^0=\cup P(\prec)_{[p]}^0$, where 
\begin{align*}
P(\prec)_{[0]}^0 & = \{\varepsilon_i-\varepsilon_j, \ \pm(\varepsilon_i+\varepsilon_j)\mid i,j\in [0],\ i\neq j\};\\
P(\prec)_{[p]}^0 & = \{\varepsilon_i-\varepsilon_j \mid\text{either } i,j\in [p],\text{ or } i,j\in [-p],\ i\neq j\} \\
& \cup \{\pm (\varepsilon_i+\varepsilon_j) \mid i\in [p],\ j\in [-p]\}
\end{align*}  
for all $[p]\neq [0]$. Hence, $\g[0]$ is of type $D$, and $\g[i]$ is of type $A$ for every $[i]\neq [0]$. Moreover,
	\[
P(\prec)^+ = \{\varepsilon_i-\varepsilon_j \mid i\prec j\} \cup \{\pm(\varepsilon_i+\varepsilon_j)\mid \pm i\prec \mp j\}.
	\]
\item If $\g=B(\infty)$, then analogously to the case of $D(\infty)$ we have $P(\prec)^0=\cup P(\prec)_{[p]}^0$, where 
\begin{align*}
P(\prec)_{[0]}^0 & = \{\pm \varepsilon_i,\  \varepsilon_i-\varepsilon_j, \ \pm(\varepsilon_i+\varepsilon_j)\mid i,j\in [0],\ i\neq j \};\\
P(\prec)_{[p]}^0 & = \{\varepsilon_i-\varepsilon_j \mid\text{either } i,j\in [p],\text{ or } i,j\in [-p],\ i\neq j\} \\
& \cup \{\pm (\varepsilon_i+\varepsilon_j) \mid i\in [p],\ j\in [-p]\}
\end{align*}  
for all $[p]\neq [0]$. Furthermore, $\g[0]$ is of type $B$, and $\g[i]$ is of type $A$ for every $[i]\neq [0]$. Moreover,
	\[
P(\prec)^+ = \{\pm \varepsilon_i\mid \pm i \prec 0\} \cup \{\varepsilon_i-\varepsilon_j \mid i\prec j\} \cup \{\pm(\varepsilon_i+\varepsilon_j)\mid \pm i\prec \mp j\}.
	\]
\item Similar statements to those of $B(\infty)$ hold for $\g=C(\infty)$ (replacing $\pm \varepsilon_i$ by $\pm 2 \varepsilon_i$).
\end{enumerate} 
\end{rem}

Throughout the paper we sometimes denote the subalgebra $\g[p]$ defined in Remark~\ref{rem:reductive.nilp.comp} by $\fsl([p])$, $\fo([p])$, $\fsp([p])$ to clearly indicate the type of $\g[p]$.

\begin{rem}[Borel subalgebras]\label{rem:linear=Borel}
We would like to point out that if we let $\prec$ be a linear order in the entire discussion above, then $\p(\prec)$ is actually a Borel subalgebra of $\g$ (see \cite[Proposition~3]{DP99}). In this case we use the notation $\fb(\prec)$ instead of $\p(\prec)$, and $\n(\prec)$ instead of $\fu(\prec)$.
\end{rem}

\subsection{Aligned modules}
Recall that all modules considered in this paper are weight modules with respect to the Cartan subalgebra $\h$ fixed in Subsection~\ref{subsec:generalities}.

\begin{defin}[$\p$-aligned modules]
Let $\p=\fl\oplus \fu$ be a parabolic subalgebra of $\g$ with reductive component $\fl$ and locally nilpotent radical $\fu$. Let $S$ be a simple weight $\fl$-module and let $\fu\cdot S=0$. The induced $\g$-module $M_\p(S):=\bU(\g)\otimes_{\bU(\p)}S$ is called a \emph{generalized Verma module}. A simple weight $\g$-module $M$ is called a \emph{$\p$-aligned module} if it is isomorphic to the unique simple quotient of $M_\p(S)$ for some simple $\fl$-module $S$.
\end{defin}

For a $\g$-module $M$ we define $M^{\fu} = \{m\in M\mid \fu\cdot m=0\}$. Moreover we let $\fu^-$ denote the opposite subalgebra of $\fu$, that is, $\fu^-$ is generated by all $\g_\alpha$ such that $\g_{-\alpha}\subseteq \fu$.

\begin{lem}\label{lem:p-ind.iff.inv.nonzero}
Let $M$ be a simple $\g$-module and suppose that $M^\fu\neq0$. Then $M^\fu$ is a simple $\fl$-module. In particular, $M$ is $\p$-aligned if and only if $M^{\fu}\neq 0$.
\end{lem}

\begin{proof}
Suppose $M^\fu\neq 0$. Since $[\p, \fu]\subseteq \fu$, it follows that $M^\fu$ is a $\p$-submodule of $M$. To see that $M^\fu$ is a simple $\fl$-module we consider the subalgebra $\bU_+(\fu^-)=\sum_{i>0}\bU^i(\fu^-)$, where $\bU^i(\fu^-)$ denotes the $i$-th step of the usual  filtration of the enveloping algebra $\bU(\fu^-)$. Then $\bU(\fu^-)=\C \oplus \bU_+(\fu^-)$, and, since $M$  is simple, we have
	\[
M=\bU(\fu^-)M^\fu=M^\fu\oplus \bU_+(\fu^-)M^\fu.
	\]
This, in particular, proves that $M^\fu$ is a simple $\fl$-module, as if $N\subseteq M^{\fu}$ is a nonzero $\fl$-submodule, then by the same argument as above we would have $M=N\oplus \bU_+(\fu^-)N$, which would imply $M^\fu\subseteq N$, since $M^\fu\cap \bU_+(\fu^-)N\subseteq M^\fu\cap \bU_+(\fu^-)M^\fu=0$. This proves the first statement.

Suppose now that $M$ is $\p$-aligned. Then there is a simple $\p$-module $S$ such that $\fu\cdot S=0$. In particular $S\subseteq M^\fu$ and hence $M^\fu\neq 0$. Conversely, if $M^\fu\neq 0$, then it is a simple $\fl$-module such that $\fu\cdot M^\fu=0$, by the first part of the lemma. Since $M$ is simple, it is clear that $\bU(\g)\otimes_{\bU(\p)} M^\fu\rightarrow M$, such that $u\otimes m \mapsto u\cdot m$, for all $u\in \bU(\fu^-)$, $m\in M^\fu$, defines a surjective homomorphism of $\g$-modules. Thus result follows.
\end{proof}

\begin{cor}\label{cor:p-ind.iff.inv.nonzero}
Let $M$ be a simple $\g$-module. Then $M$ is $\p$-aligned if and only if $M$ admits at least one $\fu$-singular weight vector.
\end{cor}
\begin{proof}
Since $M$ is assumed to be a weight module, and $M^\fu$ is an $\fl$-submodule of $M$, we have that $M^\fu$ is a weight $\fl$-module. Hence $M^\fu\neq 0$ if and only if $M^\fu_\nu\neq 0$ for some $\nu\in \h^*$. This along with Lemma~\ref{lem:p-ind.iff.inv.nonzero} implies the result.
\end{proof}

In the rest of the paper we address to the following natural question: given a simple bounded weight $\g$-module $M$ and a parabolic subalgebra $\p$, when is $M$ a $\p$-aligned module? In what follows we consider separately the cases where $M$ is integrable and non-integrable.

\section{\bf{Integrable bounded simple weight modules}}\label{sec:integrable.case}
Throughout the paper we let $V$ denote the natural representation of $\g$, which is characterized, up to isomorphism, by its support: $\Supp V = \{\varepsilon_i\}$ for $\g=\fsl(\infty)$; $\Supp V = \{0, \pm \varepsilon_i\}$ for $\g=\fo(\infty)$ and $\h=\h_B$; $\Supp V = \{\pm \varepsilon_i\}$ for $\g=\fsp(\infty)$, or $\g=\fo(\infty)$ and $\h=\h_D$ (in all cases $i\in \Z_{>0}$). The Lie algebra $\fsl(\infty)$ also admits a conatural module $V_*$, which is characterized, up to isomorphism, by its support $\Supp V_*:=-\Supp V$. Fix nonzero weight vectors $e_0\in V_0$ and $e_{\pm i}\in V_{\pm \varepsilon_i}$, consider $\g[p]$ as in Remark~\ref{rem:reductive.nilp.comp}, and set $V_{[p]}:=\langle e_i\mid i\in [p]\rangle_\C\subseteq V$. Notice that $V_{[p]}$ is isomorphic to the natural $\g[p]$-module. It was proved in \cite[\S~5]{GP18} that if $M$ is an integrable bounded simple weight $\g$-module, then $M$ is isomorphic to one of the following modules: $\Lambda_A^{\frac{\infty}{2}}V$, $S_A^\infty V$, $S_A^\infty V_*$, $S^\mu V$ or $S^\mu V_*$, if $\g=\fsl(\infty)$;  $V$, if $\g=\fsp(\infty)$; $V$, $\cS_A^B$ or $\cS_A^D$, if $\g=\fo(\infty)$. We will briefly introduce all these modules in the subsequent subsections.

\subsection{The semi-infinite fundamental modules $\Lambda_A^{\frac{\infty}{2}}V$ of $\g=\fsl(\infty)$}\label{sec:semi-inf.mod}
Let $A \subseteq \Z_{>0}$ be a semi-infinite subset, that is, $|A|=\infty$ and $|\Z_{>0}\setminus A|=\infty$. For each $n\in \Z_{>0}$, set $k_n:=A\cap [1,k]$. There is a unique, up to a constant multiplicative, embedding of $\g_n$-modules $\Lambda^{k_n}V_{n+1}\hookrightarrow \Lambda^{k_{n+1}}V_{n+2}$. We define $\Lambda_A^{\frac{\infty}{2}} V$ to be the direct limit $\g$-module $\varinjlim_{n}\Lambda^{k_n}V_{n+1}$. For two semi-infinite subsets $A, B\subseteq \Z_{>0}$, we write $A\sim_{\frac{\infty}{2}} B$ if there are finite subsets $F_A\subseteq A$, $F_B\subseteq B$ such that $A\setminus F_A = B\setminus F_B$. Hence, $\Supp \Lambda_A^{\frac{\infty}{2}} V=\{ \varepsilon_B:=\sum_{i\in B} \varepsilon_{i}\mid B\sim_{\frac{\infty}{2}} A\}$. Notice that $e_B:=\Lambda_{i\in B}e_i\in \Lambda_A^{\frac{\infty}{2}} V$ is a weight vector associated to the weight $\varepsilon_B$ (see \cite[\S~5.1]{GP18} for details).

\begin{defin}
A partial order $\prec$ on $\Z_{>0}$ is said to be compatible with a subset $A\subseteq \Z_{>0}$ if $a\in A$ and $b\in \Z_{>0}\setminus A$ implies either $a\prec b$, or $[a]=[b]$.
\end{defin}

\begin{prop}\label{prop:ext.natural.p-aligned}
Let $\prec$ be a partial order on $\Z_{>0}$. Then $\Lambda_A^{\frac{\infty}{2}} V$ is $\p(\prec)$-aligned if and only if $\prec$ is compatible with some $B\sim_{\frac{\infty}{2}} A$.
\end{prop}

\begin{proof}
Recall that $\fu(\prec)=\bigoplus_{\alpha\in P(\prec)^+}\g_\alpha$, where $P(\prec)^+=\{\varepsilon_i-\varepsilon_j\mid i\prec j\}$. It is clear that a weight $\varepsilon_B$ is not $\g_{\varepsilon_i-\varepsilon_j}$-singular if and only if $i\notin B$ and $j\in B$.

Assume that $\prec$ is compatible with some $B\sim_{\frac{\infty}{2}} A$, and let $\varepsilon_i-\varepsilon_j\in P(\prec)^+$. Since $i\prec j$, it is impossible to have $i\notin B$ and $j\in B$ (this would imply either $j\prec i$ or $[i]=[j]$). Then $\varepsilon_B$ is $\g_{\varepsilon_i-\varepsilon_j}$-singular for all $\varepsilon_i-\varepsilon_j\in P(\prec)^+$, and hence it is $\fu(\prec)$-singular. In particular, $(\Lambda_A^{\frac{\infty}{2}} V)^{\fu(\prec)}\neq 0$, and hence, by Corollary~\ref{cor:p-ind.iff.inv.nonzero}, $\Lambda_A^{\frac{\infty}{2}} V$ is $\p(\prec)$-aligned.

Conversely, suppose that $\prec$ is not compatible with any $B\sim_{\frac{\infty}{2}} A$. Then for any such $B$, there exist $j\in B$, $i\notin B$ such that $i\prec j$. This implies that $\varepsilon_i- \varepsilon_j$ is a root of $\fu(\prec)$, and by the discussed above, we know that $\varepsilon_B$ cannot be $\g_{\varepsilon_i-\varepsilon_j}$-singular. Thus no weight vector is annihilated by $\fu(\prec)$. Therefore $(\Lambda_A^{\frac{\infty}{2}} V)^{\fu(\prec)}=0$, and $\Lambda_A^{\frac{\infty}{2}} V$ cannot be $\p(\prec)$-aligned, by Corollary~\ref{cor:p-ind.iff.inv.nonzero}.
\end{proof}

\begin{prop}
Let $\prec$ be a partial order on $\Z_{>0}$. Suppose $\Lambda_A^{\frac{\infty}{2}} V$ is $\p(\prec)$-aligned and let $B\sim_{\frac{\infty}{2}} A$ be compatible with $\prec$. Then one of the following statements holds:
\begin{enumerate}
\item $\{[i]\mid [i]\cap B\neq \emptyset\}$ does not admit a maximal element, and 
	\[
(\Lambda_A^{\frac{\infty}{2}} V)^{\fu(\prec)}\cong \C
	\]
as $\fl(\prec)$-modules.
\item $[p]$ is maximal in $\{[i]\mid [i]\cap B\neq \emptyset\}$, and one of the following isomorphisms of $\fl(\prec) = \left(\bigoplus_{[i]\prec [p]} \fsl([i])\right)\oplus \fsl([p])$-modules holds: $|[p]|<\infty$, and 
	\[
(\Lambda_A^{\frac{\infty}{2}} V)^{\fu(\prec)}\cong\C\boxtimes \Lambda^{|B_{[p]}|}V_{[p]};
	\]
or $|[p]|=\infty$, and for $B_{[p]}:=B\cap [p]$ we have
	\[ 
(\Lambda_A^{\frac{\infty}{2}} V)^{\fu(\prec)} \cong \begin{cases} 
	\Lambda^{\frac{\infty}{2}}_{B_{[p]}} V_{[p]} & \text{ if } B_{[p]} \text{ is semi-infinite in } [p] \\ \\
      \Lambda^{|B_{[p]}|} V_{[p]} & \text{ if } |B_{[p]}|<\infty.
  \end{cases}
	\]
\end{enumerate}
\end{prop}

\begin{proof}
For Part (a): for any $[i]$ such that $[i]\cap B=\emptyset$ it is clear that $\fsl([i])\cdot e_B=0$. On the other hand, for any $[i]\cap B\neq \emptyset$ there is $[j]\succ [i]$ such that $[j]\cap B\neq \emptyset$. Since $B$ is compatible with $\prec$, this implies that $[i]\subseteq B$. Therefore $\fsl([i])\cdot e_B=0$, and Part (a) is proved.

For Part (b): first notice that if $[i]\prec [p]$, then $[i]\subseteq B$ since $B$ is compatible with $\prec$. In particular, $\fsl([i])\cdot e_B=0$. If there is $[i]\succ [p]$, then $[i]\cap B = \emptyset$ and $\fsl([i])\cdot e_B=0$. Finally, consider the $\fsl([p])$-module $L=\bU(\fsl([p]))\cdot e_B$. It is clear that $\varepsilon_{B_{[p]}}\in \Supp L$. Now we choose a linear order $\prec'$ on $[p]$ so that $B_{[p]}\prec [p]\setminus B_{[p]}$. By Proposition~\ref{prop:ext.natural.p-aligned}, we have that $\varepsilon_{B_{[p]}}$ is a $\fb(\prec')$-highest weight of $L$. Thus, as $\fsl([p])$-modules, we have $L\cong L_{\fb(\prec')}(\varepsilon_{B_{[p]}})$. Then either $|[p]|<\infty$ and $L\cong \Lambda^{|B_{[p]}|}V_{|[p]|}$, or $|[p]|=\infty$ and $L\cong \Lambda^{\frac{\infty}{2}}_{B_{[p]}} V_{[p]}$ if  $B_{[p]}$ is semi-infinite in $[p]$, or $L\cong \Lambda^{|B_{[p]}|} V_{[p]}$ if $|B_{[p]}|<\infty$.
\end{proof}

\begin{example}
Let $A=\{1,3,5,\ldots\}$ and consider the $\fsl(\infty)$-module $\Lambda^{\frac{\infty}{2}}_{A}V$. For the partial order $\{\ldots 5,3,1\}\prec \{2,4,6,\ldots\}$ we have $[1]\prec [2]$, where $[1]=\{\ldots 5,3,1\}$ and $[2]=\{2,4,6, \ldots\}$. Since $A$ is clearly compatible with $\prec$, we have that $\Lambda^{\frac{\infty}{2}}_{A}V$ is $\p(\prec)$-aligned, and $(\Lambda^{\frac{\infty}{2}}_{A}V)^{\fu(\prec)}$ is the $(\fsl([1])\oplus \fsl([2]))$-module generated by $e_A$. Since $A=[1]$ and $A\cap [2] = \emptyset$, we have that both $\fsl([1])$ and $\fsl([2])$ act trivially on $e_A$. Thus $(\Lambda^{\frac{\infty}{2}}_{A}V)^{\fu(\prec)}\cong \C$. This  example illustrates case (b) of previous proposition with $B = A$ and $[1] = A$ maximal in $\{[i]\mid [i]\cap B\neq \emptyset\}$.
\end{example}

\subsection{The modules $S^{\mu}V$ of $\g=\fsl(\infty)$}
Before defining the modules $S^{\mu}V$, we consider an important class of Borel subalgebras of $\g$. Recall from Remark~\ref{rem:linear=Borel} that every Borel subalgebra of $\g$ is associated to a linear order $\prec$ on $\Z_{>0}$.

\begin{defin}\label{def:Dynkin.Borel}
A Borel subalgebra $\fb(\prec)\subseteq \g$ is called a \emph{Dynkin Borel subalgebra} if and only if $(\Z_{>0},\prec)$ is isomorphic as an ordered set to $(\Z_{>0}, <)$, $(\Z_{<0}, <)$ or $(\Z, <)$. In this case we also call $\prec$ a \emph{Dynkin order}. Dynkin Borel subalgebras are the only Borel subalgebras of $\g$ for which any positive root can be written as a finite sum of simple roots.
\end{defin}

A partition of a given positive integer $\ell$ is a tuple $\mu=(\mu_1,\ldots, \mu_k)$ of positive integers such that $\mu_i\geq \mu_{i+1}$ and $\sum \mu_i=\ell$.  Let $\fb$ be the Dynkin Borel subalgebra of $\g$ associated $(\Z_{>0}, <)$, and set $\fb_n:=\fb\cap \g_n$, for every $n\in \Z_{>0}$. For a given partition $\mu$ and for $n\gg 0$, we have a unique (up to a constant multiplicative) embedding of $\g_n$-modules $L_{\fb_n}(\mu)\hookrightarrow L_{\fb_{n+1}}(\mu)$, where we identify the partition $\mu$ with the element $(\mu_1,\ldots, \mu_k, 0,\ldots, 0)\in \h_s^*$, for every $s\gg 0$. We define $S^\mu V$ to be the $\g$-module $\varinjlim_n L_{\fb_n}(\mu)$. In particular, it is not hard to see that these embeddings map $\fb_n$-highest weight vectors to $\fb_{n+1}$-highest weight vectors, and therefore $S^\mu V\cong L_{\fb}(\mu)$.

For the rest of this section we fix a partition $\mu = (\mu_1,\ldots, \mu_k)$. Moreover, for a linear order $\prec'$ on $\Z_{>0}$ which admits a left end $i_1\prec' \ldots \prec' i_k$, we define
	\[
\mu(\prec'):=\sum_{j=1}^k \mu_j\varepsilon_{i_j}.
	\]
	
For any $n\in \Z_{>0}$, we let $W_{n}$ denote the Weyl group of $\g_{n-1}$,  $\Z_n:=\{1,\ldots, n\}$, and, for a partial order $\prec$ on $\Z_{>0}$, we let $\prec_n$ denote $\prec$ restricted to $\Z_n$.

\begin{lem}\label{lem:tecnical.Smu}
Let $\prec_n$ be a partial order on $\Z_n$ and let $\prec_n'$ be a linear order on $\Z_n$. Then $\mu(\prec_n')$ is a $\fu(\prec_n)$-singular weight of $S^\mu V$ if and only if $\prec_n'$ is a refinement of $\prec_n$.
\end{lem}

\begin{proof}
Suppose that $\prec_n'$ is a refinement of $\prec_n$. Since $\mu(\prec_n')$ is clearly a $\fb(\prec_n')$-highest weight and $\fu(\prec_n)\subseteq \fb(\prec_n')$ one direction follows. Conversely, suppose $\prec_n''$ is some linear order which is a refinement of $\prec_n$. Then clearly we have $\mu(\prec_n')\in W_n\cdot \mu(\prec_n'')$. Moreover, since both weights $\mu(\prec_n')$ and $\mu(\prec_n'')$ are $\fu(\prec_n)$-singular (by assumption and by the first part) and since $(S^\mu V)^{\fu(\prec_n)}$ is a simple $\fl(\prec_n)$-module, we have $\mu(\prec_n')\in \mu(\prec_n'')+Q_{\fl(\prec_n)}$. Thus we actually can write $\mu(\prec_n')\in W_{\fl(\prec_n)}\cdot \mu(\prec_n'')$, where $W_{\fl(\prec_n)}$ denotes the Weyl group of $\fl(\prec_n)$. In other words, we have shown that the ordered set $(\Z_n, \prec_n')$ can be obtained from $(\Z_n \prec_n'')$ via an element of $W_{\fl(\prec_n)}\cong \times_{[i]} \fS_{|[i]|}$. Since $\prec_n''$ is a refinement of $\prec_n$, and any element of $\fS_{|[i]|}$ is nothing but an automorphism of the partially ordered set $(\Z_n, \prec_n)$, this proves that $\prec_n'$ is a refinement of $\prec_n$ as well.
\end{proof}

\begin{prop}\label{prop:Smu.p-aligned}
Let $\prec$ be a partial order on $\Z_{>0}$. Then $S^\mu V$ is $\p(\prec)$-aligned if and only if there are elements $i_1, \ldots, i_k\in \Z_{>0}$ (here we assume that $i_j\prec i_j$ implies $i<j$) satisfying the following property: for every $m\in \Z_{>0}\setminus\{i_1,\ldots, i_k\}$ 
\begin{equation}\label{eq:k-tail-preor.Smu}
\text{either } i_k\prec m,\text{ or } m\text{ is not comparable with } i_j,\text{ for some } j=1,\ldots, k.
\end{equation}
\end{prop}

\begin{proof}
For any $\ell\gg 0$ we define a linear order $\prec_\ell'$ on $\Z_\ell$ so that $\{i_1\prec_\ell'\cdots \prec_\ell' i_k\}\prec_\ell' (\Z_\ell \setminus \{i_1,\ldots, i_k\})$ and $\prec_\ell'$ refines $\prec_\ell$ on $(\Z_\ell\setminus \{i_1,\ldots, i_k\})$. Since $\prec$ satisfies \eqref{eq:k-tail-preor.Smu}, it is clear that $\prec_\ell'$ is a linear order which refines $\prec_\ell$ on $\Z_\ell$. By Lemma~\ref{lem:tecnical.Smu}, $\mu(\prec_\ell')$ is a $\fu(\prec_\ell)$-singular weight for any $\ell\gg 0$. Since $\mu(\prec_\ell')=\mu(\prec_{\ell'}')$ for all $\ell'\geq \ell$, we have that $\mu(\prec_\ell')$ is a $\fu(\prec)$-singular weight.

Conversely, suppose that $\prec$ does not admit elements $i_1,\ldots, i_k$ satisfying \eqref{eq:k-tail-preor.Smu}, that is, for any $k$ pairwise distinct elements $i_1,\ldots, i_k\in \Z_{>0}$ there is $i_0\in \Z_{>0}$ such that $i_0\prec i_k$ and $i_0$ is comparable with $i_j$ for all $j=1,\ldots k$. This implies the following fact: if $i_1, \ldots, i_k$ are elements satisfying \eqref{eq:k-tail-preor.Smu} in $(\Z_{\ell}, \prec)$ for some $\ell\in\Z_{> 0}$, then there is $\ell'\gg \ell$ such that the elements $i_1,\ldots, i_k$ do not satisfy \eqref{eq:k-tail-preor.Smu} in $(\Z_{\ell'},\prec)$. This implies that there is no linear order on $\Z_{\ell'}$ that refines $\prec_{\ell'}$ and which has $i_1, \ldots, i_k$ as its left end. Then, it follows from Lemma~\ref{lem:tecnical.Smu} that  $\sum_{j=1}^k \mu_j\varepsilon_{i_j}$ cannot be a $\fu(\prec_{\ell'})$-singular weight for any choice of elements $i_1,\ldots, i_k\in \Z_\ell$. In particular, no $\fu(\prec_\ell)$-singular weight is a $\fu(\prec_{\ell'})$-singular weight. Indeed, if some $\lambda^\ell$ is both $\fu(\prec_\ell)$-singular and $\fu(\prec_{\ell'})$-singular, then any $\xi\in \lambda^\ell + Q_{\fl(\ell')}$ is also $\fu(\prec_{\ell'})$-singular. But $i_1, \ldots, i_k$ are elements satisfying \eqref{eq:k-tail-preor.Smu} in $(\Z_{\ell}, \prec)$, and hence $i_1\prec \cdots \prec i_k$ is the left end of some linear order on $\Z_\ell$ that refines $\prec_n$, which implies, by Lemma~\ref{lem:tecnical.Smu}, that $\sum_{j=1}^k \mu_j\varepsilon_{i_j}$ is a $\fu(\prec_\ell)$-singular weight. Thus $\sum_{j=1}^k \mu_j\varepsilon_{i_j}\in \lambda^\ell + Q_{\fl(\ell)}\subseteq \lambda^\ell + Q_{\fl(\ell')}$, which is a contradiction.
\end{proof}

With notation as in Proposition~\ref{prop:Smu.p-aligned}, assume $i_1,\ldots, i_k\in \Z_{>0}$ are elements satisfying \eqref{eq:k-tail-preor.Smu}. Also, suppose $[i_1]\prec \cdots \prec [i_l]$ is the left end of $\{[i_j]\mid j=1,\ldots, k\}$ and write
	\[
\fl(\prec)^{ss}=\left(\bigoplus_{j=1}^l\fsl([i_j])\right)\oplus \left(\bigoplus_{[p]\succ [i_l]}\fsl([j])\right).
	\]

\begin{prop}\label{prop:Smu.inv}
With the above notation, we have the following isomorphism of $\fl(\prec)$-modules
	\[
(S^\mu V)^{\fu(\prec)}\cong \left(\boxtimes_{j=1}^l S^{\mu[i_j]}V_{[i_j]}\right)\boxtimes \C,
	\]
where $\mu[i_j]$ is the partition defined by $(\mu_l)_{l\in [i_j]}$.
\end{prop}

\begin{proof}
Notice that the existence of elements $i_1, \ldots, i_k$ satisfying \eqref{eq:k-tail-preor.Smu} is equivalent to saying the linear ordered set $(\{[i]\}, \prec)$ has a left end of size at most $k$ (in fact, the size $k$ happens when $[i_j]\neq [i_l]$ for all $j\neq l$). Take the $\fu(\prec)$-singular weight given by $\mu(\prec)=\sum_{j=1}^k \mu_{i_j}\varepsilon_j$ (see Proposition~\ref{prop:Smu.p-aligned}). Notice that $(S^\mu V)^{\fu(\prec)}$ is generated, as an $\fl(\prec)$-module, by any weight vector $v\in (S^\mu V)_{\mu(\prec)}$. In particular, $\fsl([p])$ acts trivially on $v$ for every $[p]\succ [i_k]$. On the other hand, for any $[i_j]$ we claim that the $\fsl([i_j])$-module generated by $v$ is isomorphic to $S^{\mu[i_j]}V_{[i_j]}$, where $\mu[i_j]$ is a partition defined by $(\mu_l)_{l\in [i_j]}$, and $V_{[i_j]}:=\langle e_i\mid i\in [i_j]\rangle_\C$ is the natural $\fsl([i_j])$-module. To see this, we take a linear order $\prec'$ on $\Z_{>0}$ having $i_1\prec' \cdots \prec' i_k$ as a left end, and we notice that $v$ is a $\fb(\prec')\cap \fsl([i_j])$-highest weight vector with highest weight $\sum_{i_s\in [i_j]}\mu_{i_s}\varepsilon_s$. Since $\bU(\fsl([i_j]))\cdot v$ is a simple $\fsl([i_j])$-module, the claim holds.
\end{proof}

\begin{example}
Consider a partition $\mu=(\mu_1,\mu_2,\mu_3,\mu_4,0,0,\ldots)$. Take the partial order $\{1,2,3\}\prec \{4,5,6\}\prec \{7,8,9\}\prec\cdots $. Notice that $1,2,3,4$ are elements satisfying \eqref{eq:k-tail-preor.Smu}. Thus $\sum \mu_i\varepsilon_i$ is $\fu(\prec)$-singular. The left end of $(\{[1]=[2]=[3],\ [4]\},\prec)$ is given by $[1]=\{1,2,3\}\prec [4]=\{4,5,6\}$. Hence, as an $\fl(\prec)^{ss}=\fsl([1])\oplus\fsl([4])\oplus \left(\bigoplus_{[p]\succ [4]}\fsl([j])\right)$-module, we have $(S^\mu V)^{\fu(\prec)}\cong S^{(\mu_1,\mu_2,\mu_3)}V_{[1]}\boxtimes S^{(\mu_4,0,0)}V_{[4]} \boxtimes \C$.
\end{example}

\begin{rem}
The analogs of Propositions~\ref{prop:Smu.p-aligned}~and~\ref{prop:Smu.inv} for $S^\mu V_*$ are similar to those for $S^\mu V$. Namely, ``$\prec$'' is replaced by ``$\succ$'' in Proposition~\ref{prop:Smu.p-aligned}, and ``$V_{[i_j]}$'' is replaced by ``$(V_*)_{[i_j]}$'' in Proposition~\ref{prop:Smu.inv}.
\end{rem}

\subsection{The spinor modules of $\g=\fo(\infty)$}
Set $\g_n=\fo(2n)$ and $\g_n'=\fo(2n+1)$. Recall from the theory of finite-dimensional modules that $\g_n$ admits two nonisomorphic spinor modules $\cS_n^+$ and $\cS_n^-$, and $\g_n'$ admits up to isomorphism only one spinor module $\cS_n$. If we fix the Borel subalgebras $\fb_n\subseteq \g_n$, $\fb_n'\subseteq \g_n'$ with positive roots $\{\varepsilon_i-\varepsilon_j,\ \varepsilon_k+\varepsilon_\ell,\ \varepsilon_m\mid i<j,\ k<\ell\}$, $\{\varepsilon_i-\varepsilon_j,\ \varepsilon_k+\varepsilon_\ell\mid i<j,\ k<\ell\}$, respectively, then $\cS_n^+\cong L_{\fb_n}(1/2,\ldots, 1/2)$, $\cS_n^-\cong L_{\fb_n}(1/2,\ldots, 1/2,-1/2)$, and $\cS_n\cong L_{\fb_n'}(1/2,\ldots, 1/2)$.  Thus, up to scalar, we have exactly two embeddings $\iota_{n}^\pm : \cS_{n-1}\hookrightarrow \cS_n$ of $\g_n'$-modules ($\iota^+$ maps vectors of weight $(1/2,\ldots, 1/2)\in \h_n^*$ to  vectors of weight $(1/2,\ldots, 1/2, 1/2)\in \h_{n+1}^*$ while $\iota^-$ maps vectors of weight $(1/2,\ldots, 1/2)\in \h_n^*$ to  vectors of weight $(1/2,\ldots, 1/2,-1/2)\in \h_{n+1}^*$); and a unique embedding $\cS_{n-1}^+\hookrightarrow \cS_n^+$, $\cS_{n-1}^+\hookrightarrow \cS_n^-$, $\cS_{n-1}^-\hookrightarrow \cS_n^+$, and $\cS_{n-1}^-\hookrightarrow \cS_n^-$ of $\g_n$-modules. For a subset $A\subseteq \Z_{>0}$ we define the $\h_B$-weight $\g$-module $\cS_A^B$ to be the direct limit of the $\g_n'$-modules obtained via the sequence of embeddings $\{\varphi_n: \cS_{n-1}\to \cS_n\}$, such that $\varphi_n=\iota_n^+$ if $n\in A$ and $\varphi_n=\iota_n^-$ otherwise. Similarly, we define the $\h_D$-weight  $\g$-module $\cS_A^D$ to be the direct limit of the $\g_n$-modules obtained via the sequence of embeddings $\{\varphi_n:M_{n-1}\to M_n\}$, such that $M_n=\cS_n^+$ if $n\in A$ and $M_n=\cS_n^-$ otherwise. By \cite[Theorrem~5.5]{GP18}, any simple bounded weight $\g$-module is isomorphic to $\cS_A^B$, $\cS_A^D$, or $V$. 

For every $A\subseteq \Z_{>0}$, we define $\omega_A\in \C^{\Z_{>0}}$ such that $(\omega_A)_k =\frac{1}{2}$ if $k\in A$ and $(\omega_A)_k =-\frac{1}{2}$ otherwise. For a subset  $A'\subseteq \Z_{>0}$ we write $A'\sim_B A$ if $A$ and $A'$ differ only by finitely many elements, and $A'\sim_D A$ if $A$ and $A'$ differ only by an even number of elements. It follows from \cite[\S~5.2]{GP18} that $\Supp \cS_A^B=\{\omega_{A'} \in \C^{\Z_{>0}} \mid  A'\sim_B A\}$, and $\Supp \cS_A^D=\{\omega_{A'} \in \C^{\Z_{>0}} \mid A'\sim_D A\}$.

\begin{defin}
A $\Z_2$-linear partial order $\prec$ on $\Z$ is said to be $\cS$-\emph{compatible} with $A'\subseteq \Z_{>0}$ if $\Z_{\geq 0}=A'\prec 0\prec (\Z_{>0}\setminus A')$. 
\end{defin}

It is not hard to see that $\prec$ is $\cS$-compatible with $A'\subseteq \Z_{>0}$ if and only if
\begin{equation}\label{eq:B-comp.equivalent.D-comp.}
\Z_{>0}=A'\prec (\Z_{>0}\setminus A'),\ -(\Z_{>0}\setminus A')\prec (\Z_{>0}\setminus A')\text{ and } A'\prec -A'.
\end{equation}

\begin{prop}\label{prop:spin.p-aligned}
Let $\prec$ be a $\Z_2$-linear partial order on $\Z$. Then $\cS_A^B$ (resp. $\cS_A^D$) is $\p(\prec)$-aligned if and only if $\prec$ is $\cS$-compatible with some $A'\subseteq \Z_{>0}$ for $A'\sim_B A$ (resp. $A'\sim_D A$).
\end{prop}

\begin{proof}
We will prove the result first for $\cS_A^B$. Recall that $\fu(\prec)=\bigoplus_{\alpha\in P(\prec)^+}\g_\alpha$, where 
	\[
P(\prec)^+ = \{\pm \varepsilon_i\mid \pm i \prec 0\} \cup \{\varepsilon_i-\varepsilon_j \mid i\prec j\} \cup \{\pm(\varepsilon_i+\varepsilon_j)\mid \pm i\prec \mp j\}.
	\]
Now, recall that any weight of $\cS_A^B$ is of the form $\omega_{A'}=1/2(\sum_{i\in A'} \varepsilon_{i}-\sum_{i\notin A'} \varepsilon_i)$, for some $A'\sim_B A$. Suppose that $\prec$ is $\cS$-compatible with $A'$. We claim that $\omega_{A'}$ is $\fu(\prec)$-singular. Indeed, 
\begin{enumerate}
\item If $\varepsilon_k\in P(\prec)^+$, then $k\prec 0$, which implies $k\in A'$. Thus $\omega_{A'}$ is $\g_{\varepsilon_k}$-singular.
\item If $-\varepsilon_k\in P(\prec)^+$, then $k\succ 0$, which implies $k\notin A'$. Thus $\omega_{A'}$ is $\g_{-\varepsilon_k}$-singular.
\item If $\varepsilon_i-\varepsilon_j \in P(\prec)^+$, then $i\prec j$. In particular, one cannot have $i\notin A'$ and $j\in A'$. Thus $\omega_{A'}$ is $\g_{\varepsilon_i-\varepsilon_j}$-singular.
\item If $\varepsilon_i+\varepsilon_j \in P(\prec)^+$, then $i\prec -j$. In particular, one cannot have $i,j\notin A'$. Thus $\omega_{A'}$ is $\g_{\varepsilon_i+\varepsilon_j}$-singular.
\item If $-\varepsilon_i-\varepsilon_j \in P(\prec)^+$, then $-i\prec j$. In particular, one cannot have $i,j\in A'$. Thus $\omega_{A'}$ is $\g_{-\varepsilon_i-\varepsilon_j}$-singular.
\end{enumerate}
This proves that $\omega_{A'}$ is $\fu(\prec)$-singular, and hence $\cS_A^B$ is $\p(\prec)$-aligned from Corollary~\ref{cor:p-ind.iff.inv.nonzero}. Conversely, the partial order $\prec$ not being $\cS$-compatible with any $A'\sim_B A$ implies that for every such subset $A'$ either there is $k\prec 0$ with $k\notin A'$, or there is $\ell\in A'$ such that $-\ell\prec 0$. Then either $\omega_{A'}$ is not $\g_{\varepsilon_k}$-singular, or it is not $\g_{-\varepsilon_\ell}$-singular, respectively. Thus, by Corollary~\ref{cor:p-ind.iff.inv.nonzero}, $\cS_A^B$ cannot be $\p$-aligned. 

Similarly we prove that $\cS_A^D$ is $\p(\prec)$-aligned if and only if \eqref{eq:B-comp.equivalent.D-comp.} holds for some $A'\sim_D A$. Since this is equivalent to saying $\prec$ is $\cS$-compatible with $A'$, the statement follows.
\end{proof}

Let $\prec$ be a $\Z_2$-linear partial order on $\Z$. Let $\cS_A$ be equal to $\cS_A^{B}$ or $\cS_A^{D}$, and suppose that $\p(\prec)=\fl(\prec)\oplus \fu(\prec)$ is a parabolic subalgebra for which $(\cS_A)^{\fu(\prec)}\neq 0$. By Proposition~\ref{prop:spin.p-aligned}, there is $A'\subseteq \Z_{>0}$ which is $\cS$-compatible with $\prec$ and such that $A'\sim_{B,D} A$. Set $A'_{[i]} := [i]\cap A'$, and consider
	\[
\fl(\prec)^{ss} = \left(\bigoplus_{|[i]|<\infty} \g[i]\right)\oplus \left(\bigoplus_{|[i]|=\infty} \g[i] \right).
	\]
Recall by Remark~\ref{rem:reductive.nilp.comp} that: if $\g\cong B(\infty)$, then each component $\g[i]$ is isomorphic to $A(n)$, $A(\infty)$, $B(n)$, or $B(\infty)$, with at most one component isomorphic to $B(n)$ or $B(\infty)$; if $\g\cong D(\infty)$, then each component $\g[i]$ is isomorphic to $A(n)$, $A(\infty)$, $D(n)$, or $D(\infty)$, with at most one component isomorphic to $D(n)$ or $D(\infty)$. Let $\fsl(\prec)$ denote the sum of all components of $\fl(\prec)$ which are isomorphic to either $A(n)$ or $A(\infty)$.

\begin{prop}
With the above notation, we have $\fsl(\prec)\cdot (\cS_A)^{\fu(\prec)}=0$. In particular, one of the following isomorphisms of $\fl(\prec)^{ss} = \left(\fo([0])\oplus \fsl(\prec)\right)$-modules holds:
\begin{enumerate}
\item $\g([0])\cong B(n)$, and $(\cS_A)^{\fu(\prec)}\cong \cS_n\boxtimes \C$;
\item $\g([0])\cong D(n)$, $|A'_{[0]}|\in 2\Z$ (resp. $|A'_{[0]}|\in 2\Z+1$), and $(\cS_A)^{\fu(\prec)}\cong \cS_n^+\boxtimes \C$ (resp. $(\cS_A)^{\fu(\prec)}\cong \cS_n^-\boxtimes \C$);
\item $\g([0])\cong B(\infty)$ (resp. $\g([0])\cong D(\infty)$), and $(\cS_A)^{\fu(\prec)}\cong \cS_{A'_{[0]}}^B \boxtimes \C$ (resp. $(\cS_A)^{\fu(\prec)}\cong \cS_{A'_{[0]}}^D \boxtimes \C$).
\end{enumerate}
\end{prop}

\begin{proof}
Let $[p]\neq [0]$. To see the first statement, we notice that since $A'$ is $\cS$-compatible with $\prec$, for any $i,j\in \Z_{>0}$ such that $i,j\in [p]$, we have that either $i,j\in A'$ or $i,j\in \Z_{>0}\setminus A'$. Thus $\omega_{A'}$ is  $\g_{\varepsilon_i - \varepsilon_j}$-singular. Next, suppose  $i,j\in \Z_{>0}$ are such that $i\in [p]$ and $j\in [-p]$. Since $A'$ is $\cS$-compatible with $\prec$, we may assume without loss of generality that $i\in A'$ and $j\in \Z_{>0}\setminus A'$. Hence, $\omega_{A'}$ is $\g_{\pm (\varepsilon_i + \varepsilon_j)}$-singular. Thus $\fsl([p])\cdot (\cS_A)^{\fu(\prec)}=0$, and hence $\fsl(\prec)\cdot (\cS_A)^{\fu(\prec)}=0$. To conclude the proof we look at the $\fo([0])$-module $L:=\bU(\fo([0]))\cdot (\cS_A)^{\fu(\prec)}$. We consider the case $\fo([0])\cong B(\infty)$, as the proof for the other cases is similar. It is clear that $\omega_{A'_{[0]}}\in \Supp L$, where $(\omega_{A'_{[0]}})_k=\frac{1}{2}$ if $k\in A'_{[0]}$ and $(\omega_{A'_{[0]}})_k=-\frac{1}{2}$ if $k\in [0]\setminus A'_{[0]}$. Now we choose a $\Z_2$-linear order $\prec'$ on $[0]$ so that $[0] = A'_{[0]}\prec0\prec  ([0]\setminus A'_{[0]})$. By Proposition~\ref{prop:spin.p-aligned}, we have that $\omega_{A'_{[0]}}$ is a $\fb(\prec')$-highest weight of $L$. Thus, as $\fo([0])$-modules, we have that $L\cong L_{\fb(\prec')}(\omega_{A'_{[0]}})$, which is isomorphic to $\cS_{A'_{[0]}}^B$.
\end{proof}

\subsection{The natural module of $\g=\fo(\infty)$, $\fsp(\infty)$}
Let $\prec$ be a $\Z_2$-linear partial order on $\Z$. Recall that $V$ denotes the natural $\g$-module, and that, for any $p\in \Z$, we have set $V_{[p]}:=\langle e_i\mid i\in [p]\rangle_\C\subseteq V$, which is isomorphic to the natural $\g[p]$-module.

\begin{prop}\label{prop:B,C.V.p-aligned.iff}
Let $\prec$ be a $\Z_2$-linear partial order on $\Z$. If $\g$ is $B(\infty)$, or $C(\infty)$, then $V$ is $\p(\prec)$-aligned if and only if one and only one of following statements holds:
\begin{enumerate}
\item \label{item1.prop:B,C.V.p-aligned.iff} there exists $i_0\in \Z_{>0}$ such that $i_0\prec 0$, and $i_0$ is minimal in $\Z\setminus \{0\}$.
\item \label{item2.prop:B,C.V.p-aligned.iff} there exists $j_0\in \Z_{>0}$ such that $j_0\succ 0$, and $j_0$ is maximal in $\Z\setminus \{0\}$.
\end{enumerate}
Moreover, $V^{\fu(\prec)}=V_{[i_0]}\boxtimes \C,$ as an $\fl(\prec)^{ss}=\fsl([i_0])\oplus \left( \bigoplus_{[i]\neq [i_0]}\g[i]\right)$-module if \eqref{item1.prop:B,C.V.p-aligned.iff} holds; and $V^{\fu(\prec)}=V_{[j_0]}\boxtimes \C,$ as an $\fl(\prec)^{ss}=\fsl([j_0])\oplus \left( \bigoplus_{[i]\neq [j_{0}]}\g[i]\right)$-module if \eqref{item2.prop:B,C.V.p-aligned.iff} holds.
\end{prop}

\begin{proof}
Assume $\g=C(\infty)$. Recall that $\fu(\prec)=\bigoplus_{\alpha\in P(\prec)^+}\g_\alpha$, where 
	\[
P(\prec)^+ = \{\pm 2\varepsilon_i\mid \pm i \prec 0\} \cup \{\varepsilon_i-\varepsilon_j \mid i\prec j\} \cup \{\pm(\varepsilon_i+\varepsilon_j)\mid \pm i\prec \mp j\}.
	\]
We claim that a weight of the form $\varepsilon_{i_0}$ is $\fu(\prec)$-singular if and only if $i_0$ satisfies \eqref{item1.prop:B,C.V.p-aligned.iff}. Indeed: if $j\prec i_0$ for some $j\in \Z_{>0}$, then $\varepsilon_j - \varepsilon_{i_0}\in P(\prec)^+$ and $\varepsilon_{i_0}$ is not $\g_{\varepsilon_j - \varepsilon_{i_0}}$-singular; if $-j\prec i_0$ for some $j\in \Z_{>0}$, then $-\varepsilon_j-\varepsilon_{i_0}\in P(\prec)^+$ and $\varepsilon_{i_0}$ is not $\g_{-\varepsilon_j-\varepsilon_{i_0}}$-singular; if $i_0\succ 0$, then $-2\varepsilon_{i_0}\in P(\prec)^+$ and $\varepsilon_{i_0}$ is not $\g_{-2\varepsilon_{i_0}}$-singular; finally, if $i_0\in [0]$ then $i_0$ is neither minimal nor maximal, and the result follows from the previous cases. This proves one direction of the claim. The other direction is clear. Similarly we show that a weight of the form $-\varepsilon_{j_0}$ is $\fu(\prec)$-singular if and only if $j_0$ satisfies \eqref{item2.prop:B,C.V.p-aligned.iff}. Since  $\Supp V=\{\pm\varepsilon_i\mid i\in \Z_{>0}\}$, the result follows from Corollary~\ref{cor:p-ind.iff.inv.nonzero}.

Suppose now $\g=B(\infty)$. We claim that the weight $0\in \h^*$ is not $\fu(\prec)$-singular for any partial order $\prec$. Indeed, since $[0]\neq \Z$, either there is $j\in \Z_{>0}$ such that $j\prec 0$ (which implies $0$ is not $\g_{\varepsilon_{j}}$-singular), or there is $j\in \Z_{>0}$ such that $-j\prec 0$ (which implies $0$ is not $\g_{-\varepsilon_{j}}$-singular). This implies that the only possible $\fu(\prec)$-singular weights of $V$ are of the form $\pm\varepsilon_i$, for some $i\in \Z_{>0}$. Now the result follows similarly to the case $C(\infty)$.

Finally, it is easy to see that if \eqref{item1.prop:B,C.V.p-aligned.iff} holds then  \eqref{item2.prop:B,C.V.p-aligned.iff} does not, and vice versa. The last statement follows from the fact that $[i_0]\neq [0]$ and $[j_0]\neq [0]$, and  $V^{\fu(\prec)}$ is generated, as an $\fl(\prec)^{ss}$-module, by either $e_{i_0}$ or $e_{j_0}$, respectively.
\end{proof}

\begin{prop}\label{prop:D.V.p-aligned.iff}
Let $\prec$ be a $\Z_2$-linear partial order on $\Z$. If $\g=D(\infty)$, then $V$ is $\p(\prec)$-aligned if and only if one and only one of following statements holds:
\begin{enumerate}
\item \label{item1.prop:D.V.p-aligned.iff} there exist $i_0\in \Z_{>0}$ such that $i_0$ is minimal in $\Z\setminus \{0\}$.
\item \label{item2.prop:D.V.p-aligned.iff} there exist $j_0\in \Z_{>0}$ such that  $j_0$ is maximal in $\Z\setminus \{0\}$.
\end{enumerate}
Moreover, $V^{\fu(\prec)}=V_{[i_0]}\boxtimes \C,$ as an $\fl(\prec)^{ss}=\fsl([i_0])\oplus \left( \bigoplus_{[i]\neq [i_0]}\g[i]\right)$-module if \eqref{item1.prop:B,C.V.p-aligned.iff} holds; and $V^{\fu(\prec)}=V_{[j_0]}\boxtimes \C,$ as an $\fl(\prec)^{ss}=\fsl([j_0])\oplus \left( \bigoplus_{[i]\neq [j_{0}]}\g[i]\right)$-module if \eqref{item2.prop:B,C.V.p-aligned.iff} holds.
\end{prop}

\begin{proof}
Similarly to the proof of Proposition~\ref{prop:B,C.V.p-aligned.iff}, we prove that $\varepsilon_{i_0}$ (resp. $-\varepsilon_{i_0}$) is $\fu(\prec)$-singular if and only if $i_0$ is minimal (resp. maximal) in $\Z\setminus \{0\}$. The moreover part is also similar to Proposition~\ref{prop:B,C.V.p-aligned.iff}.
\end{proof}

\section{\bf{Non-integrable bounded simple weight modules}}\label{sec:non-integrable.case}
It follows from \cite{GP18} that $\g$ admits non-integrable bounded simple weight modules if and only if $\g=\fsl(\infty),\ \fsp(\infty)$. In this section we study the structure of such modules. Let $\g=\varinjlim_{n} \g_n$ be one of these Lie algebras. By \cite[Corollary~4.3]{GP18}), if $M$ is a bounded simple weight module, then $M$ is locally simple, that is, $M=\varinjlim_n M_n$, where $M_n$ is a simple $\g_n$-module for every $n\in \Z_{>0}$. Thus any integrable bounded simple weight module $M$ is a direct limit of simple finite-dimensional $\g_n$-modules, and for any splitting Borel subalgebra $\fb=\varinjlim_{n} \fb_n$ (resp. parabolic subalgebra $\p=\varinjlim_{n} \p_n$), $M$ is a direct limit of simple $\fb_n$-highest weight (resp. $\p_n$-aligned) $\g_n$-modules. Recall from the previous sections that it may be the case that $M$ itself is not a $\fb$-highest weight (resp. $\p$-aligned) module. When this is the case we say $M$ is a \emph{pseudo $\fb$-highest weight (resp. pseudo $\p$-aligned) module}. In particular, any integrable bounded simple weight module is either a $\fb$-highest weight (resp. $\p$-aligned) module or a pseudo $\fb$-highest weight (resp. pseudo $\p$-aligned) module. On the other hand, if $M$ is non-integrable, it may be the case that $M$ is neither $\fb$-highest weight ($\p$-aligned) nor pseudo $\fb$-highest weight (resp. pseudo $\p$-aligned) as we shall see in this section.

Most of the results of this section regarding (pseudo) $\fb$-highest weight modules are equivalent to results that have been established in \cite[\S~6]{GP18}. Since our terminology is slightly different from that of \cite{GP18}, we add their proofs here for reader's convenience.

\subsection{The modules $X_\fsl(\mu)$ of $\g=\fsl(\infty)$}\label{sec:X_sl}

For details on what follows we refer to \cite[\S~6]{GP18}. Recall that the Weyl algebra of differential operators of $\C[x_1,x_2\ldots]$ is given by $\cD_\infty=\C[x_1,x_2\ldots,\partial_1,\partial_2\ldots]$ with relations $\partial_i x_j - x_j\partial_i = \delta_{i,j}$, $x_ix_j=x_jx_i$, $\partial_i\partial_j=\partial_j\partial_i$. For any $\mu\in \C^{\Z_{>0}}$, define $\bx^\mu=\Pi_{i\geq 1}x_i^{\mu_i}$, and consider the space
	\[
\cF_{\fsl}(\mu) = \{\bx^\mu p \mid p\in \C[x_1^{\pm 1}, x_2^{\pm 1},\ldots],\ \deg p=0\} = \Span\{\bx^\lambda\mid \lambda-\mu\in Q_\g\},
	\]
where we are identify the root lattice $Q_\g$ with the set $\{(\lambda_1,\lambda_2,\ldots)\in \C_{\rm f}^{\Z_{>0}}\mid \sum \lambda_i=0\}$, where $\C_{\rm f}^{\Z_{>0}}$ is the set of all sequences of complex numbers admitting only finitely many nonzero entries.  The space $\cF_{\fsl}(\mu) $ is a bounded weight $\g$-submodule, where the action is defined through the homomorphism of associative algebras $\bU(\g)\rightarrow \cD_{\infty}$ given by $\g_{\varepsilon_i - \varepsilon_j}\ni E_{ij}\mapsto x_i\partial_j$, where $E_{ij}$ is the elementary matrix with $1$ in the $i,j$-position and zeros elsewhere. Let $\Int(\mu)$, $\Int^+(\mu)$ and $\Int^-(\mu)$ denote the subsets of $\Z_{> 0}$ composed by all $i$ such that $\mu_i\in \Z$, $\mu_i\in \Z_{\geq 0}$ and $\mu_i\in \Z_{<0}$, respectively. Then
	\[
V_{\fsl}(\mu)=\Span\{\bx^\lambda\mid \lambda-\mu\in Q_\g,\ \Int^+(\mu)\subseteq \Int^+(\lambda)\}
	\]
is a $\g$-submodule of $\cF_\fsl(\mu)$. Moreover, we have that $V_\fsl(\mu')\subseteq V_\fsl(\mu)$ if and only if $\mu-\mu'\in Q_\g$ and $\Int^+(\mu)\subseteq \Int^+(\mu')$ (in fact, if $\ell\in\Int^+(\mu')\setminus \Int^+(\mu)$, then we cannot get $x_\ell^{\mu_\ell}$ from $x_\ell^{\mu_\ell'}$ by acting with $x_\ell^n\partial_\ell^k$, $n,k\in \Z_{\geq 0}$). Now we set $V_\fsl(\mu)^+:=0$ if $\Int^-(\mu)=\emptyset$, and
	\[
V_\fsl(\mu)^+:=\sum_{V_\fsl(\mu')\subsetneq V_\fsl(\mu)} V_\fsl(\mu')
	\]
otherwise (the summation above runs over all $\mu'\in \C^{\Z_{>0}}$ such that $\mu-\mu'\in Q_\g$ and $\Int^+(\mu)\subsetneq \Int^+(\mu')$). Finally, we define
	\[
X_\fsl(\mu) := V_\fsl(\mu)/V_\fsl(\mu)^+.
	\]

\begin{rem}\label{rem:properties.X_sl} The following properties can be found in \cite[\S~6.2.1]{GP18}:
\begin{enumerate}
\item Any non-integrable bounded simple weight $\g$-module is isomorphic to $X_{\fsl}(\mu)$ for some $\mu\in \C^{\Z_{>0}}$;
\item If ${\bf x}^\lambda\in X_\fsl(\mu)$, then $(x_i\partial_j)\cdot {\bf x}^\lambda=0$ if and only if either $\lambda_i=-1$ or $\lambda_j=0$. Moreover, ${\bf x}^\lambda$ is a weight vector of weight $\lambda\in \h^*$;
\item $\Supp X_\fsl(\mu)=\{\lambda\in \C^{\Z_{>0}}\mid \lambda \sim_{\fsl} \mu\}$, where $\lambda\sim_{\fsl}\mu$ means that $\lambda-\mu\in Q_\g$, and $\Int^+(\mu)=\Int^+(\lambda)$. In particular, if $\lambda\sim_{\fsl}\mu$ then $\Int^-(\mu) = \Int^-(\lambda)$;
\item $X_{\fsl}(\mu)\cong X_{\fsl}(\mu')$ if and only if $\mu \sim_{\fsl} \mu'$, or $\{\mu,\mu'\}=\{{\bf 0}, -{\bf 1}\}$;
\item For any $\mu^n\in \C^n$ we can define a simple $\g_{n-1}$-module  $X_{\fsl}(\mu^n)$ similarly to $X_\fsl(\mu)$.
\end{enumerate}
\end{rem}

Let $\fb$ be a Borel subalgebra of $\g$ and $\p$ be a parabolic subalgebra of $\g$. We now examine whether $X_\fsl(\mu)$ is a $\fb$-highest (resp. pseudo $\fb$-highest) weight module or a $\p$-aligned (resp. pseudo $\p$-aligned) module. 

Recall that Borel subalgebras of $\g$ are in bijection with linear orders on $\Z_{>0}$ (see Remark~\ref{rem:linear=Borel}). Given a linear order $\prec$ on $\Z_{>0}$, $i_0\in\Z_{>0}$ and $c\in \C$, we define $\varepsilon(\prec, i',c)\in \C^{\Z_{>0}}$ such that $\varepsilon(\prec, i',c)_i=-1$ if $i\prec i'$, $\varepsilon(\prec, i',c)_i=0$ if $i'\prec i$ and $\varepsilon(\prec, i',c)_{i'}=c$. Also, given a semi-infinite subset $J \subseteq \Z_{>0}$ which is compatible with $\prec$ (i.e. $i\in J$ and $j\notin J$ implies $i\prec j$), we define $\varepsilon(\prec, J)\in \C^{\Z_{>0}}$ such that $\varepsilon(\prec, J)_{i}=-\delta_{i\in J}$, where $\delta_{i\in J}=1$ if $i\in J$ and $0$ if $i\notin J$. In what follows, for any $c\in\C$ we let ${\bf c}$ denote the constant tuple $(c,c,c,\ldots)$ that may be finite or infinite (it will be clear from the context which case we are considering).

\begin{prop}\label{prop:X(mu).hw.iff}
Suppose that $X_{\fsl}(\mu)$ is not a trivial module, and let $\prec$ be a linear order on $\Z_{>0}$. Then $X_{\fsl}(\mu)$ is a $\fb(\prec)$-highest weight module if and only if either $\mu\sim_{\fsl} \varepsilon(\prec, i',c)$ for some $i'\in \Z_{>0}$ and $c\in \C$, or $\mu\sim_{\fsl} \varepsilon(\prec, J)$ for some semi-infinite subset $J\subseteq \Z_{>0}$.
\end{prop}

\begin{proof}
It is clear that the $\g$-module $X_{\fsl}(\varepsilon(\prec, i',c))$ is isomorphic to $L_{\fb(\prec)}(\varepsilon(\prec, i',c))$, since ${\bf x}^{\varepsilon(\prec, i',c)}$ is a $\fb(\prec)$-highest weight vector of $X_{\fsl}(\varepsilon(\prec, i',c))$. Since $X_{\fsl}(\mu)\cong X_{\fsl}(\varepsilon(\prec, i',c))$ if and only if $\varepsilon(\prec, i',c)\in \Supp X_{\fsl}(\mu)$ (see Remark~\ref{rem:properties.X_sl}), we have proved one direction. (Notice that the same argument also works for $X_{\fsl}(\varepsilon(\prec, J))$).

Conversely, we start by noting that $X_\fsl(-{\bf 1})\cong X_\fsl({\bf 0})\cong \C$. Then we may assume that $\mu\neq -{\bf 0}$ and $-{\bf 1}$. Let $\lambda\in \Supp X_{\fsl}(\mu)$ be such that ${\bf x}^\lambda$ is annihilated by $\n(\prec)$. If $\lambda_{i'}\notin \{-1,0\}$ for some $i'\in \Z_{>0}$, then it is easy to see that $\lambda=\varepsilon(\prec, i',\lambda_{i'})$, and we are done. Suppose now $\lambda_i\in \{-1,0\}$ for every $i\in \Z_{>0}$. Let $J=\{i\in \Z_{>0}\mid \lambda_i=-1\}$. If $|J|<\infty$ (resp. $|\Z\setminus J|<\infty$), then we set $i' :=\max \{i\in J\}$ (resp. $i'=\min\{i\in\Z_{>0}\setminus J\}$), and it is easy to see that $\lambda=\varepsilon(\prec, i', -1)$ (resp. $\lambda=\varepsilon(\prec, i', 0)$). Finally, if $J$ is semi-infinite, then $\lambda=\varepsilon(\prec,J)$. This concludes the proof.
\end{proof}

\begin{example}\label{ex:varepsilon(prec, J)}
Consider the linear order $1\prec 3\prec 5\prec\cdots\prec 6\prec 4\prec 2$ on $\Z_{>0}$, and let $\mu=(-1,0,-1,0,\ldots)\in \C^{\Z_{>0}}$. It is clear that $X_{\fsl}(\mu)$ is not trivial and that $\mu$ is a $\fb(\prec)$-highest weight. However, one can see that $\mu\nsim \varepsilon(\prec, i', c)$ for any choice of $i'\in \Z_{>0}$, $c\in \C$, but clearly $\mu\sim_{\fsl} \varepsilon(\prec, J)$, where $J=\{2k+1\mid k\in \Z_{\geq 0}\}$. The next result shows that $X_\fsl(\mu)$ dos not admit such highest weights with respect to a Dynkin Borel subalgebra (recall Definition~\ref{def:Dynkin.Borel}).
\end{example}

\begin{cor}
Suppose that $X_{\fsl}(\mu)$ is not a trivial module, and let $\prec^d$ be a Dynkin linear order on $\Z_{>0}$. Then $X_{\fsl}(\mu)$ is a $\fb(\prec^d)$-highest weight module if and only if $\mu\sim_{\fsl} \varepsilon(\prec^d, i',c)$ for some $i'\in \Z_{>0}$ and $c\in \C$.
\end{cor}

\begin{proof}
Regarding the proof of Proposition~\ref{prop:X(mu).hw.iff}, it only remains to consider the case where $J=\{i\in \Z_{>0}\mid \lambda_i=-1\}$ is semi-infinite. Since we are considering a Dynkin order, there must exist a pair of elements $i',j'\in \Z_{>0}$ such that $\lambda_{i'}=0$ and $\lambda_{j'}=-1$, where $j'\prec^d i'$ and there is no element between $j'$ and $i'$. We claim that $\lambda_j=0$ for all $i'\prec^d j$ (resp. $\lambda_{i}=-1$ for all $i\prec^d j'$). Indeed, if $j\in \Z_{>0}$ is such that $i'\prec^d j$, then $(x_{i'}\partial_j)\cdot {\bf x}^\lambda=0$ implies $\lambda_j=0$. Similarly, if $i\in \Z_{>0}$ is such that $i \prec^d j'$, then $(x_i\partial_{j'})\cdot {\bf x}^\lambda=0$ implies $\lambda_i=-1$. Thus $\lambda = \varepsilon(\prec^d, i',0)$. In particular, we have shown that $\varepsilon(\prec^d, i',c)\in \Supp X_{\fsl}(\mu)$, for some $i'\in \Z_{>0}$ and some $c\in \C$, and thus the result follows.
\end{proof}

\begin{cor}\label{cor:X(mu).hw.iff}
Suppose that $X_{\fsl}(\mu)$ is not a trivial module, and let $\prec$ be a linear order on $\Z_{>0}$. If $X_{\fsl}(\mu)$ is a $\fb(\prec)$-highest weight module, then $\mu\sim_{\fsl} \varepsilon(\prec^d, i',c)$ for some $i'\in \Z_{>0}$, $c\in \C$, and some Dynkin order $\prec^d$ on $\Z_{>0}$. In particular, $X_{\fsl}(\mu)\cong L_{\fb(\prec^d)}(\varepsilon(\prec^d, i',c))$.
\end{cor}
\begin{proof}
Assume that $X_{\fsl}(\mu)\cong L_{\fb(\prec)}(\lambda)$ for some $\lambda\in \h^*$. Since $\dim X_{\fsl}(\mu)_\lambda<\infty$ for every $\lambda\in \Supp X_{\fsl}(\mu)$, it follows from \cite[Proposition~5-(iii)]{DP99} that there is some Dynkin Borel subalgebra $\fb(\prec^d)$ for which $L_{\fb(\prec)}(\lambda)\cong L_{\fb(\prec^d)}(\lambda)$. By Proposition~\ref{prop:X(mu).hw.iff}, this is the case if and only if $\lambda=\varepsilon(\prec^d, i',c)$ for some $i'\in \Z_{>0}$ and $c\in \C$. Thus the result is proved.
\end{proof}

\begin{example}
In Example~\ref{ex:varepsilon(prec, J)} we had $X_\fsl(\mu)\cong  L_{\fb(\prec)}(\varepsilon(\prec, J))$. But notice that we also have $X_{\fsl}(\mu)\cong L_{\fb(\prec^d)}(\varepsilon(\prec^d, 1,-1))$, where $\cdots\prec^d 3\prec^d 1\prec^d 2\prec^d 4\prec \cdots$ is a Dynkin order on $\Z_{>0}$.
\end{example}

For any $\mu\in \C^{\Z_{>0}}$ and any $J\subseteq \Z_{>0}$ we set $\mu_{J} := (\mu_\ell)_{\ell\in J}\in \C^{|J|}$. Moreover, given two sequences $\gamma=(\gamma_i)\in \C^{\Z_{>0}}$ and $\eta=(\eta_i)\in \C^{\Z_{>0}}$, we write $T(\gamma)=T(\eta)$ if $\gamma_i=\eta_i$ for all but finitely many indices $i\in \Z_{>0}$. Finally, for convenience we define 
	\[
I=\{i\in \Z_{>0}\mid i\notin \Int(\mu)\},\quad F_- = \Int^-(\mu),\quad F_+ = \Int^+(\mu).
	\]
In particular, $\Z_{>0}=F_-\cup I\cup F_+$.

\begin{cor}
Let $\prec$ be a linear order on $\Z_{>0}$. Then $X_{\fsl}(\mu)$ is a pseudo $\fb(\prec)$-highest weight if and only if $|I|\leq 1$ and either $T(\mu_{F_-})\neq T(-{\bf 1})$ or $T(\mu_{F_+})\neq T({\bf 0})$.
\end{cor}

\begin{proof}
If $|I|\geq 2$, then $X_\fsl(\mu)$ is not locally $\fb(\prec)$-highest weight. On the other hand, if $|I|\leq 1$, then $X_\fsl(\mu)$ is clearly locally $\fb(\prec)$-highest weight. Moreover, by Proposition~\ref{prop:X(mu).hw.iff}, $X_\fsl(\mu)$ is $\fb(\prec)$-highest weight if and only if $\mu\sim_{\fsl} \varepsilon(\prec, i',c)$ for some $i'\in \Z_{>0}$ and $c\in \C$, or $\mu\sim_{\fsl} \varepsilon(\prec, J)$ for some semi-infinite subset $J\subseteq \Z_{>0}$. But this is the case if and only if $T(\mu_{F_-})= T(-{\bf 1})$ and $T(\mu_{F_+}) = T({\bf 0})$. Hence the result follows.
\end{proof}

\begin{example}
Take the linear order $1\prec 3\prec 5\prec\cdots\prec 6\prec 4\prec 2$ on $\Z_{>0}$, and consider $\mu=(-1,0,c,0,1,1,\ldots)\in \C^{\Z_{>0}}$ for some $c\in \C$. Since $T(\mu_{F_+})\neq T({\bf 0})$, it follows that $X_\fsl (-1,0,c,0,1,1,\ldots)$ is does not have a $\fb(\prec)$-highest weight, but is a pseudo $\fb(\prec)$-highest weight $\g$-module.
\end{example}

Let $\prec$ be a partial order on $\Z_{>0}$, and consider the parabolic subalgebra $\p(\prec)=\fl(\prec)\oplus \fu(\prec)$. Recall from Remark~\ref{rem:reductive.nilp.comp} that
	\[
\fl(\prec)^{ss} = \bigoplus \fsl([i]).
	\]
Next we define the subset $S(\prec)\subseteq \C^{\Z_{>0}}$ such that $\varepsilon\in S(\prec)$ if and only if $[i]\prec [j]$ implies either $\varepsilon_{[i]}=-{\bf 1}$ or $\varepsilon_{[j]}={\bf 0}$. Also, for any $\mu\in \C^{\Z_{>0}}$, we set $S_\mu(\prec):=\{\lambda\in S(\prec)\mid \lambda\sim_{\fsl}\mu\}$.

\begin{prop}\label{prop:X_sl.p-top.}
Let $\prec$ be a partial order on $\Z_{>0}$. Then $X_{\fsl}(\mu)$ is $\p(\prec)$-aligned if and only if  $\mu\sim_{\fsl} \varepsilon\in S_\mu(\prec)$. Moreover, one of the following statements holds:
\begin{enumerate}
\item \label{item1:prop:X_sl.p-top.} $\varepsilon_{[i']}\notin \{{\bf 0}, -{\bf 1}\}$ for a unique $[i']$, and either $|[i']|\geq 2$ and
	\[
X_{\fsl}(\mu)^{\fu(\prec)} \cong X_{\fsl}(\varepsilon_{[i']})\boxtimes \C
	\]
as $\left(\fsl([i'])\oplus \left(\bigoplus_{[j]\neq [i]} \fsl([j])\right)\right)$-modules, or $|[i']|=1$ and $X_{\fsl}(\mu)^{\fu(\prec)} \cong \C$ as $\fl(\prec)^{ss}$-modules. 
\item\label{item2:prop:X_sl.p-top.}  $\varepsilon_{[i]}\in \{{\bf 0}, -{\bf 1}\}$ for all $[i]$, and $X_{\fsl}(\mu)^{\fu(\prec)}\cong \C$ as $\fl(\prec)^{ss}$-modules.
\end{enumerate}
\end{prop}

\begin{proof}
It is clear that $\varepsilon\in \Supp X_{\fsl}(\mu)$ is $\fu(\prec)$-singular if and only if $\varepsilon\in S_\mu(\prec)$. Then it follows from Corollary~\ref{cor:p-ind.iff.inv.nonzero} that $X_{\fsl}(\mu)$ is $\p(\prec)$-aligned if and only if  $\mu\sim_{\fsl} \varepsilon\in S_\mu(\prec)$. Now we have the following possibilities
\begin{enumerate}
\item $\varepsilon_{[i']}\notin \{{\bf 0}, -{\bf 1}\}$ for some $[i']$. In this case, we must have that $\varepsilon_{[j]}= -{\bf 1}$ for all $[j]\prec [i']$, and $\varepsilon_{[j]}= {\bf 0}$ for all $[i']\prec [j]$. In particular, such $[i']$ is unique.
\item $\varepsilon_{[i]}\in \{{\bf 0}, -{\bf 1}\}$ for all $[i]$. In this case, $[i]\prec [j]$ implies either $\varepsilon_{[i]}= {\bf 0}$ and $\varepsilon_{[j]}= {\bf 0}$; or  $\varepsilon_{[i]}= -{\bf 1}$ and $\varepsilon_{[j]}= -{\bf 1}$; or $\varepsilon_{[i]}= -{\bf 1}$ and $\varepsilon_{[j]}= {\bf 0}$.
\end{enumerate}	

Notice that $\fl(\prec)^{ss}$ acts trivially on ${\bf x}^\varepsilon$ if $\varepsilon_{[i]}\in \{{\bf 0}, -{\bf 1}\}$ for all $[i]$. Thus, in such a case we have an isomorphism $X_{\fsl}(\mu)^{\fu(\prec)}\cong \C$ as $\fl(\prec)^{ss}$-modules. On the other hand, if $\varepsilon_{[i']}\notin \{{\bf 0}, -{\bf 1}\}$ for some $[i']$, then $\fsl([j])$ acts trivially on ${\bf x}^\varepsilon$ for every $[j]\neq [i']$. Hence, if $|[i']|=1$, then it is clear that $X_\fsl(\mu)^{\fu(\prec)}\cong \C$ as $\fl(\prec)^{ss}$-modules. If $|[i']|\geq 2$, then $\bU(\fsl([i']))\cdot {\bf x}^\varepsilon$ generates an $\fsl([i'])$-module isomorphic to $X_{\fsl}(\varepsilon_{[i]})$. Thus, in this case we have that $X_{\fsl}(\mu)^{\fu(\prec)}\cong X_{\fsl}(\varepsilon_{[i']})\boxtimes \C$ as $\left(\fsl([i'])\oplus \left(\bigoplus_{[j]\neq [i']} \fsl([j])\right)\right)$-modules.
\end{proof}

\begin{cor}\label{cor:X_sl.p-aligned.iff}
Let $\prec$ be a partial order on $\Z_{>0}$. Then $X_{\fsl}(\mu)$ is $\p(\prec)$-aligned if and only if one of the following statements holds:
\begin{enumerate}
\item  There is a unique $[i']$ such that $\Z_{>0}=(F_-\setminus [i'])\prec [i']\prec (F_+\setminus [i'])$, and $|F_-\setminus [i']|=\infty$ implies $T(\mu_{F_-\setminus [i']})=T(-{\bf 1})$, and $|F_+\setminus [i']|=\infty$ implies $T(\mu_{F_+\setminus [i']})=T({\bf 0})$.
\item $\Z_{>0}=F_-\prec F_+$, and $|F_-|=\infty$ implies $T(\mu_{F_-})=T(-{\bf 1})$, and $|F_+|=\infty$ implies $T(\mu_{F_+})=T({\bf 0})$.
\end{enumerate}
\end{cor}

Let $\mu\in \C^{\Z_{>0}}$. For any subset $Z\subseteq \Z_{>0}$ such that $\Int^-(\mu_Z)=\mu_{Z}$ we set
\begin{equation}\label{eq:s^pm}
s^-(\mu_Z):=\sum_{i\in Z}(\mu_i+1).
\end{equation}
Similarly, if $\Int^+(\mu_Z)=\mu_{Z}$, then we set 
	\[
s^+(\mu_Z):=\sum_{i\in Z} \mu_i. 
	\]


\subsection*{The shadow of $X_{\fsl}(\mu)$}

Let $M$ be a simple weight $\g$-module and let $\lambda\in \Supp M$ be any fixed weight. Let $\alpha\in \D$ and consider the set $\sm_\alpha^\lambda=\{q\in \R\mid \lambda+q\alpha\in \Supp M\}\subseteq \R$. Notice that $\sm_\alpha^\lambda$ may be of four types: bounded from both, above and below, bounded from above but not from below, bounded from below but not from above, and unbounded from both, above and below. This induces a decomposition of $\D$ into four subsets:
\begin{align*}
\D_\lambda^\cF & = \{\alpha\in \D\mid \sm_\alpha^\lambda \text{ is bounded in both directions}\} \\
\D_\lambda^\cI & = \{\alpha\in \D\mid \sm_\alpha^\lambda \text{ is unbounded in both directions}\} \\
\D_\lambda^+ & = \{\alpha\in \D\mid \sm_\alpha^\lambda \text{ is bounded from above but not from below}\} \\
\D_\lambda^- & = \{\alpha\in \D\mid \sm_\alpha^\lambda \text{ is bounded from below but not from above}\}.
\end{align*}
It is easy to see that $\alpha\in \D_\lambda^+$ if and only if $-\alpha\in \D_\lambda^-$. In other words, we have $\D_\lambda^- = -\D_\lambda^+$. 

It was proved in \cite[Theorem~2]{DP99} that this decomposition of $\D$ does not depend on $\lambda$, but only on $M$. Therefore we write $\D_M^\bullet$ instead of $\D_\lambda^\bullet$. Consider the induced decomposition of $\g$
	\[
\g=\left(\g_M^\cF + \g_M^\cI\right)\oplus \g_M^+\oplus \g_M^-,
	\]
where $\g_M^\cF:=\h\oplus \left(\bigoplus_{\alpha\in \D_M^\cF}\g_\alpha\right)$, $\g_M^\cI:=\h\oplus \left(\bigoplus_{\alpha\in \D_M^\cI}\g_\alpha\right)$, and $\g_M^\pm:=\bigoplus_{\alpha\in \D_M^\pm}\g_\alpha$. The triple $(\g_M^\cI, \g_M^+, \g_M^-)$ is called the \emph{shadow} of $M$ on $\g$.  If $\g$ is a finite dimensional reductive Lie algebra, then it is true that $\p_M=(\g_M^\cF+\g_M^\cI)\oplus \g_M^+$ is a parabolic subalgebra of $\g$, its reductive component is given by $\g_M^\cF+\g_M^\cI$, its nilradical is $\g_M^+$, and there exist a natural surjection
	\[
\bU(\g)\otimes_{\bU(\p)} M^{\g_M^+}\rightarrow M,
	\]
where $M^{\g_M^+}$ is the simple $(\g_M^\cF+\g_M^\cI)$-submodule of $M$ spanned by all elements that are annihilated by $\g_M^+$. Moreover, $\Supp M = \Supp \bU(\g)\otimes_{\bU(\p_M)} M^{\g_M^+}$. This is known as the Fernando-Futorny parabolic induction theorem. On the other hand, if $\g=\fsl(\infty)$, $\fo(\infty)$, $\fsp(\infty)$, then $\p_M$ still defines a parabolic subalgebra of $\g$, however, the Fernando-Futorny parabolic induction theorem may not hold in general (since the nilpotent subalgebra $\g_M^+$ is infinite-dimensional, it may be the case that $M^{\g_M^+}=0$ (see \cite[Example~5]{DP99})). 

For any $\mu\in \C^{\Z_{>0}}$,  we have that the reductive (resp. locally nilpotent) component of $\p_{X_\fsl(\mu)}$ is given by $\fl_{X_\fsl(\mu)} :=\g_{X_{\fsl}(\mu)}^\cF + \g_{X_{\fsl}(\mu)}^\cI,\quad (\text{resp. } \fu_{X_\fsl(\mu)} :=\g_{X_{\fsl}(\mu)}^+)$.

\begin{lem}\label{lem:shadow.X(mu)}
The shadow of $X_{\fsl}(\mu)$ is given by
\begin{align*}
\D_{X_\fsl(\mu)}^\cF & = \{\varepsilon_i-\varepsilon_j\mid i,j \in F_+,\text{ or } i,j \in F_-\}; \\
\D_{X_\fsl(\mu)}^+ & = \{\varepsilon_i-\varepsilon_j\mid i \in F_- \text{ and }  j \notin F_-,\text{ or } j \in F_+ \text{ and } i \notin F_+ \}; \\
\D_{X_\fsl(\mu)}^- & = -\D_{X_\fsl(\mu)}^+; \\
\D_{X_\fsl(\mu)}^\cI & = \{\varepsilon_i-\varepsilon_j \mid i,j \notin \Int (\mu)\}.
\end{align*}
\end{lem}

\begin{proof}
This follows from the following facts: $x_i\partial_j$ acts locally finite on $X_{\fsl}(\mu)$ if and only if $i\in F_-$ or $j\in F_+$; $x_i\partial_j$ acts injectively on $X_{\fsl}(\mu)$ if and only if $i\notin F_-$ and $j\notin F_+$.
\end{proof}

A subset $F\subseteq \D$ is called \emph{indecomposable} if for every pair of roots $\alpha,\beta\in F$ there exists a sequence $\alpha = \alpha_1,\ \alpha_2,\ldots, \alpha_n = \beta$ in $F$ such that $\alpha_{i}+\alpha_{i+1}\in F$ for all $i=1,\ldots, n-1$. A subset which is not indecomposable is called \emph{decomposable}.

\begin{lem}\label{lem:levi.p_X_sl}
Then the following statements hold:
\begin{enumerate}
\item  If $\D_{X_\fsl(\mu)}^\cI\neq \emptyset$, then $\D_{X_\fsl(\mu)}^\cI$ is indecomposable. In particular, $[i]=I$ for every $i\in I$.
\item  If $\D_{X_\fsl(\mu)}^\cF\neq \emptyset$, then $\D_{X_\fsl(\mu)}^\cF=(\D_{X_\fsl(\mu)}^\cF)^-\sqcup (\D_{X_\fsl(\mu)}^\cF)^+$, where $(\D_{X_\fsl(\mu)}^\cF)^\pm=\{\varepsilon_i-\varepsilon_j\mid i,j\in F_\pm \}$ are both indecomposable. In particular $[i]=F_\pm$, for every $i\in F_\pm$, respectively.
\end{enumerate}
\end{lem}

\begin{proof}
Part (a): the set $\D_{X_\fsl(\mu)}^\cI$ is clearly additively closed. To see that it is indecomposable, note that for any pair of roots $(\varepsilon_i-\varepsilon_j), (\varepsilon_k-\varepsilon_\ell)\in \D_{X_\fsl(\mu)}^\cI$, the sequence $\alpha_1=(\varepsilon_i-\varepsilon_j),\ \alpha_2=(\varepsilon_k-\varepsilon_i),\ \alpha_3=(\varepsilon_j-\varepsilon_k),\ \alpha_4=(\varepsilon_k-\varepsilon_\ell)$ lies in $\D_{X_\fsl(\mu)}^\cI$ and is such that $\alpha_i+\alpha_{i+1}\in \D_{X_\fsl(\mu)}^\cI$. Part (b) follows similarly.
\end{proof}

Let $\prec_\mu$ denote the partial order on $\Z_{>0}$ corresponding to the parabolic subalgebra $\p_{X_\fsl(\mu)}$. By Lemma~\ref{lem:levi.p_X_sl}, we have that $|[i]|>1$ for at most three different classes, which coincide with the sets $I$ and $F_\pm$. Let $\fsl(I)$, $\fsl(F_\pm)$ denote the corresponding subalgebras. In particular, we have that $\fl_{X_\fsl(\mu)}^{ss}=\fsl(I)\oplus \fsl(F_-)\oplus \fsl(F_+)$. Moreover, since the set of roots of $\fu_{X_\fsl(\mu)}$ is $\D_{X_\fsl(\mu)}^+$, we also have that $F_-\prec_\mu I\prec_\mu F_+$.  The next result answers the following natural question: when does the Fernando-Futorny parabolic induction theorem hold for $X_{\fsl}(\mu)$? Or equivalently, when is $X_{\fsl}(\mu)$ a $\p_{X_\fsl(\mu)}$-aligned module?

\begin{prop}
One of the following statements holds:
\begin{enumerate}
\item $|F_\pm|<\infty$, and $X_{\fsl}(\mu)$ is $\p_{X_\fsl(\mu)}$-aligned. In this case, there is $\varepsilon_{I}\in \C^{|I|}$ such that
	\[
X_{\fsl}(\mu)^{\fu_{X_\fsl(\mu)}}\cong X_{\fsl}(\varepsilon_{I})\boxtimes \C
	\]
as $\fsl(I)\oplus (\fsl(F_-)\oplus \fsl(F_+))$-modules.
\item $I\neq \emptyset$, and $X_{\fsl}(\mu)$ is $\p_{X_\fsl(\mu)}$-aligned if and only if $|F_-|=\infty$ implies $T(\mu_{F_-})=T(-{\bf 1})$, and $|F_+|=\infty$ implies $T(\mu_{F_+})=T({\bf 0})$. In this case, there is $\varepsilon_{I}\in \C^{|I|}$ such that
	\[
X_{\fsl}(\mu)^{\fu_{X_\fsl(\mu)}}\cong X_{\fsl}(\varepsilon_{I})\boxtimes \C
	\]
as $\fsl(I)\oplus (\fsl(F_-)\oplus \fsl(F_+))$-modules; 
\item $I=\emptyset$, and $X_{\fsl}(\mu)$ is $\p_{X_\fsl(\mu)}$-aligned if and only if $|F_-|=\infty$ implies $T(\mu_{F_-})=T(-{\bf 1})$, and $|F_+|=\infty$ implies $T(\mu_{F_+})=T({\bf 0})$. In this case, either $s^-(\mu_{F_-})+s^+(\mu_{F_+})=0$ and $X_{\fsl}(\mu)^{\fu_{X_\fsl(\mu)}}\cong \C$, or $s^-(\mu_{F_-})+s^+(\mu_{F_+})\neq 0$ and there is $\varepsilon_{F_-}\in \Z^{|F_-|}$ (resp. $\varepsilon_{F_+}\in \Z^{|F_+|}$) such that
	\[	
X_{\fsl}(\mu)^{\fu_{X_\fsl(\mu)}} \cong X_{\fsl}(\varepsilon_{F_-})\boxtimes \C\quad (\text{resp. } X_{\fsl}(\mu)^{\fu_{X_\fsl(\mu)}} \cong \C\boxtimes X_{\fsl}(\varepsilon_{F_+}))
	\]
as $\fsl(F_-)\oplus \fsl(F_+)$-modules.
\end{enumerate}
\end{prop}

\begin{proof}
This follows from Proposition~\ref{prop:X_sl.p-top.}.
\end{proof}

Let $\mu\in \C^{\Z_{>0}}$ and $\mu^n=(\mu_1,\ldots, \mu_n)\in \C^n$. A partial order $\prec_n$ on $\Z_n$ is said to be $\fsl$-\emph{compatible} with $\mu^n$ if $\prec_n$ coincides with $\prec_{\mu^n}$ on $\Z_n\setminus \Int(\mu^n)$, and $\Int^-(\mu^n)\prec_n (\Z_n\setminus \Int(\mu^n))\prec_n \Int^+(\mu^n)$. A partial order $\prec$ on $\Z_{>0}$ is said to be \emph{locally $\fsl$-compatible} with $\mu$ if $\prec_n$ is $\fsl$-compatible with $\mu^n$ for every $n\in \Z_{>0}$.

\begin{rem}\label{rem:nice.prop}
The subalgebra $\p_{X_\fsl(\mu)}$ satisfies the following nice property: $X_{\fsl}(\mu^n)$ is $\p_{X_\fsl(\mu^n)}$-aligned for every $n\in \Z_{>0}$. In particular, $X_{\fsl}(\mu)$ is pseudo $\p_{X_\fsl(\mu)}$-aligned if and only if it is not $\p_{X_\fsl(\mu)}$-aligned. The next result shows that the parabolic subalgebra $\p_{X_\fsl(\mu)}$ plays a similar role in the theory of aligned modules to that played by Dynkin Borel subalgebras in the theory of highest weight modules.
\end{rem}

\begin{prop}\label{prop:X_sl.pseudo.p.aligned}
Let $\prec$ be a partial order on $\Z_{>0}$. Then $X_{\fsl}(\mu)$ is pseudo $\p(\prec)$-aligned if and only if $\prec$ is locally $\fsl$-compatible with $\mu$ and $X_{\fsl}(\mu)$ is not $\p(\prec)$-aligned.
\end{prop}

\begin{proof}
Notice first that $X_{\fsl}(\mu)\cong \varinjlim_n X_{\fsl}(\mu^n)$, and $X_{\fsl}(\mu^n)$ is $\p_{X_\fsl(\mu^n)}$-aligned for every $n\in \Z_{>0}$.  We claim that $X_{\fsl}(\mu^n)$ is $\p(\prec_n)$-aligned if and only if $\prec_n$ is $\fsl$-compatible with $\mu^n$. Indeed, if $I\neq \emptyset$, then for $n\gg 0$ we have that $\mu^n\sim_{\fsl} \varepsilon$, where $\varepsilon_i=-1$ for all $i\in \Int^-(\mu^n)$, and $\varepsilon_i=0$ for all $i\in \Int^+(\mu^n)$. If $\prec_n$ is $\fsl$-compatible with $\mu^n$, then $i\prec_n j$ implies either $i,j\in \Int^-(\mu^n)$, or $i,j\in \Int^+(\mu^n)$, or $i\in \Int^-(\mu^n)$ and $j\in I\cup \Int^+(\mu^n)$, or $i\in I$ and $j\in \Int^+(\mu^n)$. Then, for all cases, either $\varepsilon_i=-1$ or $\varepsilon_j=0$ and $\varepsilon$ is $\fu(\prec_n)$-singular. Now this direction follows from Corollary~\ref{cor:p-ind.iff.inv.nonzero}. On the other hand, if $\prec_n$ is not $\fsl$-compatible with $\mu^n$, then either there are $i,j\in I$ with $i\prec_n j$, which clearly gives a contradiction, since for any weight $\lambda\in \Supp X_\fsl(\mu^n)$, we have that $\lambda_i,\lambda_j\notin \Z$, and hence there is no weight which is $\g_{\varepsilon_i-\varepsilon_j}$-singular; or there exist $i\in \Int^-(\mu^n)$, $j\in I\cup \Int^+(\mu^n)$ with $j\prec i$, or $i\in I$, $j\in \Int^+(\mu^n)$ with $j\prec i$, and in both cases there is no weight $\lambda$ of $X_{\fsl}(\mu^n)$ which is $\g_{\varepsilon_j-\varepsilon_i}$-singular (notice that we are using the following fact: if $\nu\in \Supp X_{\fsl}(\eta)$, then $\Int^\pm(\nu) = \Int^\pm (\eta)$). Thus, for the case $I\neq \emptyset$, we have shown that $X_\fsl(\mu^n)$ is $\p(\prec_n)$-aligned if and only if $\prec_n$ is $\fsl$-compatible with $\mu^n$. Finally, if $I=\emptyset$, then $\mu^n\sim_{\fsl} \varepsilon(\prec_{\mu^n}, i',c)$, for some $i'\in \Z_{>0}$ and some $c\in \Z$. Hence, similarly to the previous case, we show that $X_\fsl(\mu^n)$ is $\p(\prec_n)$-aligned if and only if $\prec_n$ is $\fsl$-compatible with $\mu^n$. Thus $X_\fsl(\mu)$ is locally $\p(\prec)$-aligned if and only if $\prec$ is locally $\fsl$-compatible with $\mu$, and the statement follows.
\end{proof}

\subsection{The integrable modules $S_A^\infty V$ and $S_A^\infty V_*$ of $\g=\fsl(\infty)$}\label{sec:S_infty.case}
In the previous section we have answered question \eqref{Q2} (see Introduction) for $\g$-modules of the form $X_{\fsl}(\mu)$. In particular, this gives an answer to this question for all non-integrable simple bounded weight $\g$-modules. Notice however that $X_{\fsl}(\mu)$ is not necessarily non-integrable. The goal of this section is precisely to cover a case where $X_{\fsl}(\mu)$ is integrable. 

We discuss now the only class of integrable bounded simple weight modules that was not considered in Section~\ref{sec:integrable.case}. Namely, the modules $S_A^\infty V$,  $S_A^\infty V_*$. To define such modules we let $A=\{a_1,a_2,\ldots\mid a_i \leq a_{i+1}\}\subseteq \Z_{>0}$. For each $a_n\in \Z_{>0}$, there is a unique, up to a constant multiplicative, embedding of $\g_{n-1}$-modules $S^{a_n}V_{n}\hookrightarrow S^{a_{n+1}}V_{n+1}$. We define the $\g$-module $S_A^{\infty} V$ to be the direct limit $\varinjlim_{n}S^{a_n}V_{n}$. The module $S_A^{\infty} V_*$ is defined similarly.

The next lemma along with results of Sections~\ref{sec:integrable.case}~and~\ref{sec:X_sl} and Remark~\ref{rem:X.iso.Symm_*} below answer question \eqref{Q2} for all integrable simple bounded weight $\g$-modules.

\begin{lem}\label{lem:X.iso.Symm}
There is an isomorphism of $\g$-modules 
	\[
S_A^\infty V\cong X_{\fsl}(a_1,a_2-a_1,a_3-a_2,\ldots).
	\]
\end{lem}

\begin{proof}
As $\g_{n-1}$-modules, we have $X_{\fsl}(a_1,a_2-a_1,\ldots, a_{n}-a_{n-1})\cong X_{\fsl}(a_n,0,\ldots, 0)\cong S^{a_n} V_n$, for every $n\in \Z_{>0}$. Since $X_{\fsl}(a_1,a_2-a_1,a_3-a_2,\ldots)\cong \varinjlim_n X_{\fsl}(a_1,a_2-a_1,\ldots, a_{n}-a_{n-1})$, the result follows.
\end{proof}

The next proposition is a well known fact which follows from results of \cite{DP99} (see also \cite[Proposition~5.2]{GP18}). We include its proof here for completeness.

\begin{prop}\label{prop:Symm.hw}
With the above notation, assume that the sequence $(a_1,a_2,\ldots)$ does not stabilize. Then $S_A^\infty V$ is not a highest weight module with respect to any Borel subalgebra.
\end{prop}

\begin{proof}
Using Proposition~\ref{prop:X(mu).hw.iff} and the fact that $\Z_{>0}=\Int^+(a_1,a_2-a_1,\ldots)$, we see that $S_A^\infty V$ is $\fb(\prec)$ highest weight if and only if $(a_1,a_2-a_1,\ldots)\sim_{\fsl} \varepsilon(\prec, i', c)$ for some $i'\in \Z_{>0}$ and $c\in \Z$. Since $\Z_{>0}=\Int^+(a_1,a_2-a_1,\ldots)$, we see that this is the case if and only if $i'$ is the minimal element of $\Z_{>0}$ and $T(a_1,a_2-a_1,\ldots)=T({\bf 0})$. But the latter equality is impossible, since by hypothesis we always can find $i\gg 0$ such that $a_i<a_{i+1}$.
\end{proof}

\begin{example}[Parabolic subalgebra with two infinite blocks]
Let $\prec$ be the partial order on $\Z_{>0}$ such that $[1]:=\{\cdots, 5, 3, 1\}\prec [2]:=\{2, 4, 6,\cdots\}$, and let $A = \{1,1,2,2, 3,3,4,4,\ldots\}\subseteq \Z_{\geq 0}$. Then $S_A^\infty V\cong X_{\fsl}(1,0,1,0,1,0\ldots)$. Since $(1,0,1,0,1,0\ldots)_{[2]}={\bf 0}$, $S_A^\infty V$ is clearly $\p(\prec)$-aligned. Moreover, as $\fsl([1])\oplus \fsl([2])$-modules, we have 
	\[
(S_A V)^{\fu(\prec)}\cong X_{\fsl}(1,1,\ldots)\boxtimes \C\cong S_B^\infty V\boxtimes \C,
	\]
where $B=\{1,2,3,4,\ldots\}$.
\end{example}

\begin{rem}\label{rem:X.iso.Symm_*}
Similarly to Lemma~\ref{lem:X.iso.Symm}, one can prove that 
	\[
S_A^\infty V_*\cong X_{\fsl}(- {\bf 1} - (a_1,a_2-a_1, a_3-a_2,\ldots)).
	\]
In particular, the analog of  Proposition~\ref{prop:Symm.hw} for $S_A^\infty V_*$ implies that such a module is not a highest weight module with respect to any Borel subalgebra of $\g$.
\end{rem}

\subsection{The modules $X_\fsp(\mu)$ of $\g=\fsp(\infty)$}

In this section we follow the notation used in Section~\ref{sec:X_sl}.  For any $\mu\in \C^{\Z_{>0}}$ we let $\g$ act on 
	\[
\cF_{\fsp}(\mu) = \{\bx^\mu p \mid p\in \C[x_1^{\pm 1}, x_2^{\pm 1},\ldots],\ \deg p\in 2\Z\} = \{{\bf x}^\lambda\mid \lambda-\mu\in Q_\g\}
	\]
via the homomorphism of associative algebras $\bU(\g)\rightarrow \cD_{\infty}$, defined by: $g_{\varepsilon_i-\varepsilon_j}\mapsto x_i\partial_j$, $g_{\varepsilon_i+\varepsilon_j}\mapsto x_ix_j$, $g_{2\varepsilon_i}\mapsto 1/2x_i^2$, $g_{-\varepsilon_i-\varepsilon_j}\mapsto -\partial_i\partial_j$, $g_{-2\varepsilon_i}\mapsto -(1/2)\partial_i^2$, where $g_\alpha\in \g_{\alpha}$ are appropriate nonzero vectors. Now the definition of the $\g$-module $X_{\fsp}(\mu)$ is similar to that of $X_{\fsl}(\mu)$ (see Section~\ref{sec:X_sl}). 

\begin{rem}\label{rem:properties.X_sp} The following properties can be found in \cite[\S~6.2.2]{GP18}:
\begin{enumerate}
\item Any non-integrable bounded simple weight $\g$-module is isomorphic to $X_{\fsp}(\mu)$ for some $\mu\in \C^{\Z_{>0}}$;
\item The vector ${\bf x}^\lambda$ is a weight vector of weight $\lambda + {\bf \frac{1}{2}}\in \h^*$;
\item $\Supp X_\fsp(\mu)=\{\lambda + {\bf \frac{1}{2}}\mid \lambda \sim_{\fsp} \mu\}$, where $\lambda\sim_{\fsp}\mu$ means that $\lambda-\mu\in Q_\g$ and $\Int^+(\mu)=\Int^+(\lambda)$. In particular, if $\lambda\sim_{\fsp}\mu$ then $\Int^-(\mu) = \Int^-(\lambda)$;
\item $X_{\fsp}(\mu)\cong X_{\fsp}(\mu')$ if and only if $\mu'\sim_{\fsp} \mu$;
\item For any $\mu^n\in \C^n$ we can define a simple $\g_n$-module  $X_{\fsp}(\mu^n)$ similarly to $X_\fsp(\mu)$.
\end{enumerate}
\end{rem}

Recall that Borel subalgebras of $\g$ are in bijection with $\Z_2$-linear orders on $\Z$ (see Remark~\ref{rem:linear=Borel}). Let $\prec$ be a $\Z_2$-linear order on $\Z$ and define $\omega(\prec)\in \C^{\Z_{>0}}$ such that $\omega(\prec)_i=-1$ if $i\prec 0$, and $\omega(\prec)_i=0$ if $i\succ 0$. Moreover, recall that $\varepsilon_i\in\C^{\Z_{>0}}$ is being identified with the sequence that has $1$ at the i-th entry and zeros elsewhere.

\begin{prop}\label{prop:X_sp(mu).hw.iff}
Let $\prec$ be a $\Z_2$-linear order on $\Z$. Then $X_{\fsp}(\mu)$ is a $\fb(\prec)$-highest weight module if and only if either $\mu\sim_{\fsp} \omega(\prec)$ or $\mu\sim_{\fsp} \omega(\prec)\mp \varepsilon_{i'}$, where the latter holds only if $i'$ (resp. $-i'$) is a maximal element in $\{\pm i\mid \pm i\prec 0\}$.
\end{prop}

\begin{proof}
If ${\bf x}^\lambda\in X_{\fsp}(\mu)$ is a $\fb(\prec)$-highest weight vector, then it is clear that $\lambda_i\in \{-2,-1\}$ (resp. $\lambda_i\in \{0,1\}$) for all $i\in \Z_{>0}$ such that $i\prec 0$ (resp. $i\succ 0$). Assume that $i\prec 0$ is not a maximal element in $\{\pm i\mid \pm i\prec 0\}$. Then there exists $j\in \Z_{>0}$ such that $\g_{\varepsilon_i\pm \varepsilon_j}\subseteq \n(\prec)$. This implies that $\lambda_i=-1$. Similarly, if $-i\prec 0$ in not maximal in $\{\pm i\mid \pm i\prec 0\}$, then we prove that $\lambda_i=0$. In particular, $\lambda=\omega(\prec)$ in this case. Finally, it is clear that if $i'$ (resp. $-i'$) is a maximal element in $\{\pm i\mid \pm i\prec 0\}$, then $\lambda$ may be of the form $\omega(\prec)\mp \varepsilon_{i'}$, respectively. The other implication is obvious.
\end{proof}

Recall the definition of $s^\pm(\mu_{F_\pm})$ in \eqref{eq:s^pm}, and, for each $n\in \Z_{>0}$, set 
	\[
s_n^+(\mu_{F_+}):=\sum_{i=1}^n \mu_i, \quad s_n^{-}(\mu_{F_-}):=\sum_{i=1}^n(\mu_i+1).
	\]
	
\begin{cor}\label{cor:X_sp.h.w.I=empty}
Let $\prec$ be a $\Z_2$-linear order on $\Z$. Then $X_{\fsp}(\mu)$ is $\fb(\prec)$-highest weight if and only if $I=\emptyset$ and $s^+(\mu_{F_+}) + s^-(\mu_{F_-})\in \{0,\pm 1\}$.
\end{cor}

\begin{cor}
Let $\prec$ be a $\Z_2$-linear order on $\Z$. Then $X_{\fsp}(\mu)$ is a pseudo $\fb(\prec)$-highest weight $\g$-module if and only if $I=\emptyset$ and the sequence $(s_n^+(\mu_{F_+}) + s_n^-(\mu_{F_-}))_{n\in \Z_{>0}}\in \{0,\pm 1\}^{\Z_{>0}}$ does not converge.
\end{cor}
\begin{proof}
If $I\neq \emptyset$, then $X_\fsl(\mu)$ is not locally $\fb(\prec)$-highest weight. On the other hand, it is clear that $X_\fsp(\mu)$ is locally $\fb(\prec)$-highest weight if and only if  $I=\emptyset$ and $(s_n^+(\mu_{F_+}) + s_n^-(\mu_{F_-}))_{n\in \Z_{>0}}\in \{0,\pm 1\}^{\Z_{>0}}$. Then the statement follows from Corollary~\ref{cor:X_sp.h.w.I=empty}.
\end{proof}

\subsection*{$\p$-aligned analysis}
Let $\prec$ be a $\Z_2$-linear partial order on $\Z$, and consider the parabolic subalgebra $\p(\prec)=\fl(\prec)\oplus \fu(\prec)$. Recall from Remark~\ref{rem:reductive.nilp.comp} that 
\begin{align*}
P(\prec)^+ & = \{2\varepsilon_i\mid i \prec 0\}\cup \{-2\varepsilon_i\mid -i \prec 0\} \cup \{\varepsilon_i-\varepsilon_j \mid i\prec j\} \\
& \cup \{\varepsilon_i+\varepsilon_j\mid i\prec -j\}\cup \{-\varepsilon_i-\varepsilon_j\mid - i\prec j\},
\end{align*}
and that $P(\prec)^0=\cup P(\prec)_{[p]}^0$, where 
\begin{align*}
P(\prec)_{[0]}^0 & = \{\pm \varepsilon_i,\  \varepsilon_i-\varepsilon_j, \ \pm(\varepsilon_i+\varepsilon_j)\mid i,j\in [0],\ i\neq j \}; \\
P(\prec)_{[p]}^0 & = \{\varepsilon_i-\varepsilon_j \mid\text{either } i,j\in [p],\text{ or } i,j\in [-p],\ i\neq j\} \\
& \cup \{\pm(\varepsilon_i+\varepsilon_j)\mid i\in [p],\ j\in [-p]\},
\end{align*}  
for all $[p]\neq [0]$. Define $R(\prec)\subseteq \C^{\Z_{>0}}$ so that $\omega\in R(\prec)$ if and only if $\omega_{[i]\cap\Z_{>0}}=-{\bf 1}$ for all $[i]\prec [0]$, and $\omega_{[i]\cap\Z_{>0}}={\bf 0}$ for all $[i]\succ [0]$. Moreover, we set $R_\mu(\prec):=\{\omega\in R(\prec)\mid \omega\sim_{\fsp} \mu\}$. Now assume that $|[0]|>1$. If the set $\{[i]\prec [0]\}$ admits a maximal element, then we define $\omega(+)\in R(\prec)$ such that $\omega(+)_{[0]\cap\Z_{>0}}={\bf 0}$. Similarly, if $\{[i]\succ [0]\}$ admits a minimal element, then we define $\omega(-)\in R(\prec)$ such that $\omega(-)_{[0]\cap\Z_{>0}}=-{\bf 1}$.

\begin{prop}\label{prop:X_sp(mu).p-aligned}
Let $\prec$ be a $\Z_2$-linear partial order on $\Z$. Then $X_{\fsp}(\mu)$ is $\p(\prec)$-aligned if and only if one of the following statements holds:
\begin{enumerate}
\item $\{[i]\prec [0]\}$ does not admit a maximal element, $\{[i]\succ [0]\}$ does not admit a minimal element, and $\mu\sim_{\fsp} \omega\in R_\mu(\prec)$.
\item $\{[i]\prec [0]\}$ admits a maximal element $[i']$, $\{[i]\succ [0]\}$ does not admit a minimal element, and either $\mu\sim_{\fsp} \omega(+) -\varepsilon_{i_0}$ for  some $i_0\in [i']$, or $\mu\sim_{\fsp} \omega\in R_\mu(\prec)$.
\item $\{[i]\succ [0]\}$ admits a minimal element $[i']$, $\{[i]\prec [0]\}$ does not admit a maximal element, and either $\mu\sim_{\fsp} \omega(-) +\varepsilon_{i_0}$ for  some $i_0\in [i']$, or $\mu\sim_{\fsp} \omega\in R_\mu(\prec)$.
\item $\{[i]\prec [0]\}$ admits a maximal element $[i_1]$, $\{[i]\succ [0]\}$ admits a minimal element $[i_2]$, and either $\mu\sim_{\fsp} \omega(+) -\varepsilon_{i_0}$ for  some $i_0\in [i_1]$, or $\mu\sim_{\fsp} \omega(-) +\varepsilon_{i_0}$ for  some $i_0\in [i_2]$, or $\mu\sim_{\fsp} \omega\in R_\mu(\prec)$.
\end{enumerate}
\end{prop}

\begin{proof}
If ${\bf x}^\lambda\in X_{\fsp}(\mu)$ is a $\fu(\prec)$-singular weight vector, then it is clear that $\lambda_i\in \{-2,-1\}$ (resp. $\lambda_i\in \{0,1\}$) for all $i\in \Z_{>0}$ such that $[i]\prec [0]$ (resp. $[i]\succ [0]$). If there exists $[i']\prec [0]$ for which $\lambda_{i'}=-2$, then it is easy to see that $\lambda_i=-1$ (resp. $\lambda_i=0$) for all $[i]\prec [i']$ (resp. $[i]\succ [i']$). In particular, this implies that $[i']$ is maximal in $\{[i]\prec [0]\}$. Similarly, if there is $[i']\succ [0]$ for which $\lambda_{i'}=1$, then $\lambda_i=-1$ (resp. $\lambda_i=0$) for all $[i]\prec [i']$ (resp. $[i]\succ [i']$), and $[i']$ must be minimal in $\{[i]\succ [0]\}$. This along with Corollary~\ref{cor:p-ind.iff.inv.nonzero} shows Part~(a). On the other hand, if $[i']$ is maximal in $\{[i]\prec [0]\}$, then either $\lambda_{i_0}=-2$ for a unique $i_0\in [i']$ (to see this we use the fact that $\lambda$ must be $\g_{\varepsilon_i+\varepsilon_j}$-singular, for all $i,j \in [i']$) and $\lambda=\omega(+)- \varepsilon_{i_0}$, or $\lambda_i=-1$ for all $i\in [i']$ and $\lambda\sim_{\fsp} \omega\in R_\mu(\prec)$, which concludes Part~(b) (again by Corollary~\ref{cor:p-ind.iff.inv.nonzero}). Part~(c) and (d) are similar, and the converse is obvious.
\end{proof}

Let $\prec$ be a $\Z_2$-linear partial order on $\Z$, assume $\p(\prec)=\fl(\prec)\oplus \fu(\prec)$ is a parabolic subalgebra for which $X_{\fsp}(\mu)$ is $\p(\prec)$-aligned, and consider the decomposition (see Remark~\ref{rem:reductive.nilp.comp})
	\[
\fl(\prec)^{ss} := \fsp([0])\oplus \bigoplus_{[i]\neq [0]} \fsl([i]).
	\]
	
\begin{cor}\label{cor:X(sp).top}
With the above notation, one of the following statements holds:
\begin{enumerate}
\item $\mu\sim_{\fsp} \omega\in R_\mu(\prec)$, and
	\[
X_{\fsp}(\mu)^{\fu(\prec)} \cong X_{\fsp}(\omega_{[0]})\boxtimes \C,
	\]
as $\fsp([0])\oplus \left(\bigoplus_{[i]\neq [0]} \fsl([i])\right)$-modules.
\item $\mu\sim_{\fsp} \omega(+) - \varepsilon_{i_0}$ for some $[i_0]\prec [0]$, and 
	\[
X_{\fsp}(\mu)^{\fu(\prec)} \cong \C\boxtimes X_{\fsl}(-{\bf 1}-\varepsilon_{i_0})\boxtimes \C,
	\]
as $\fsp([0])\oplus\fsl([i_0]) \oplus \left(\bigoplus_{[i]\neq [0],[i_0]} \fsl([i])\right)$-modules.
\item $\mu\sim_{\fsp} \omega(-) + \varepsilon_{i_0}$ for some $[i_0]\succ [0]$, and 
	\[
X_{\fsp}(\mu)^{\fu(\prec)} \cong \C\boxtimes X_{\fsl}(\varepsilon_{i_0})\boxtimes \C,
	\]
as $\fsp([0])\oplus\fsl([i_0]) \oplus \left(\bigoplus_{[i]\neq [0],[i_0]} \fsl([i])\right)$-modules.
\item $\mu\sim_{\fsp} \omega(+)$, or $\mu\sim_{\fsp} \omega(-)$, and 
	\[
X_{\fsp}(\mu)^{\fu(\prec)} \cong \C,
	\]
as $\fl(\prec)^{ss}$-modules.
\end{enumerate}
\end{cor}

\subsection*{The shadow of $X_{\fsp}(\mu)$}
Recall the definition of the shadow of a simple weight module given in Section~\ref{sec:X_sl}. For any $\mu\in \C^{\Z_{>0}}$ we have that the reductive component (resp. locally nilpotent) component of $\p_{X_{\fsp}(\mu)}$ is given by 
	\[
\fl_{X_\fsp(\mu)} := \g_{X_{\fsp}(\mu)}^\cF + \g_{X_{\fsp}(\mu)}^\cI\quad (\text{resp.} \fu_{X_\fsp(\mu)} := \g_{X_{\fsp}(\mu)}^+).
	\]

\begin{lem}\label{lem:shadow.X_sp(mu)}
The shadow of $X_{\fsp}(\mu)$ is given by
\begin{align*}
\D_{X_\fsp(\mu)}^\cF & = \{\varepsilon_i-\varepsilon_j\mid i,j \in F_+,\text{ or } i,j \in F_-\}\cup \{\pm(\varepsilon_i+\varepsilon_j) \mid i\in F_-,\ j\in F_+\};  \\
\D_{X_\fsp(\mu)}^+ & = \{2\varepsilon_i\mid i\in F_-\}\cup\{-2\varepsilon_i\mid i\in F_+\}\cup \{\varepsilon_i+\varepsilon_j \mid i \in F_-, \ j\in F_-\cup I\} \\
&\cup  \{-\varepsilon_i-\varepsilon_j \mid i \in F_+, \ j\in I\cup F_+\} \\
& \cup \{\varepsilon_i-\varepsilon_j\mid i \in F_-, \ j \in I\cup F_+,\text{ or } i \in I,\ j \in F_+ \}; \\
\D_{X_\fsp(\mu)}^- & = -\D_{X_\fsp(\mu)}^+; \\
\D_{X_\fsp(\mu)}^\cI & = \{\pm 2\varepsilon_i\mid i\in I\} \cup \{\varepsilon_i-\varepsilon_j \mid \{i,j\}\cap \Int(\mu)=\emptyset \} \\
& \cup \{\pm(\varepsilon_i+\varepsilon_j) \mid \{i,j\} \cap \Int(\mu) = \emptyset  \}.
\end{align*}
\end{lem}

\begin{proof}
Similar to Lemma~\ref{lem:shadow.X(mu)}.
\end{proof}

\begin{lem}\label{lem:levi.p_X}
The following statements hold:
\begin{enumerate}
\item  If $\D_{X_\fsp(\mu)}^\cI\neq \emptyset$, then $\D_{X_\fsp(\mu)}^\cI$ is indecomposable.
\item  If $\D_{X_\fsp(\mu)}^\cF\neq \emptyset$, then $\D_{X_\fsp(\mu)}^\cF$ is indecomposable.
\end{enumerate}
\end{lem}

\begin{proof}
Part (a) is obvious. Part (b): using similar arguments to those of Lemma~\ref{lem:levi.p_X}, we see that $\D_X^\cF$ is indecomposable. Moreover, if $\{\pm (\varepsilon_i+\varepsilon_j),\ \pm (\varepsilon_i-\varepsilon_j)\}\subseteq \D_{X_\fsp(\mu)}^\cF$, then $\pm 2\varepsilon_i\in \D_{X_\fsp(\mu)}^\cF$ which is a contradiction. Hence the subalgebra associated to $\D_{X_\fsp(\mu)}^\cF$  is isomorphic to $\fsl(F)$.
\end{proof}

Let $\prec_\mu$ denote the $\Z_2$-linear partial order on $\Z$ corresponding to the parabolic subalgebra $\p_{X_\fsp(\mu)}$. By Lemma~\ref{lem:levi.p_X}, we have that $|[i]|>1$ for at most three different classes, which we denote by $[0]=\{\pm i\mid i\in \Z_{>0}\text{ and } i\in I\}\cup \{0\}$, $[0_-]=\{ i\mid i \in F_-\}\cup \{ -i\mid i \in F_+\}$ and $[0_+]=\{ i\mid i \in F_+\}\cup \{ -i\mid i \in F_-\}$. Let $\fsp(I)$ denote the subalgebra associated to $[0]$ and $\fsl(F)$ denote the subalgebra corresponding to $[0_-] \cup [0_+]$. In particular, we have that $\fl_{X_\fsp(\mu)}^{ss}=\fsp(I)\oplus \fsl(F)$. Moreover, since the set of roots of $\fu_{X_\fsp(\mu)}$ is $\D_{X_\fsp(\mu)}^+$, it is easy to see that $[0_-]\prec_\mu[0]\prec_\mu [0_+]$.  The following result is the Fernando-Futorny parabolic induction theorem for $X_{\fsp}(\mu)$.

\begin{prop}
The module $X_{\fsp}(\mu)$ is $\p_{X_\fsp(\mu)}$-aligned if and only if one of the following statements holds:
\begin{enumerate}
\item $I\neq \emptyset$, $\mu\sim_{\fsp} \omega\in R_\mu(\prec_\mu)$, and
	\[
X_{\fsp}(\mu)^{\fu_{X_\fsp(\mu)}} \cong X_{\fsp}(\omega_{I})\boxtimes \C
	\]
as $\fsp(I)\oplus \fsl(F)$-modules. In this case, $\omega$ can be chosen so that for any fixed $i'\in I$: $\omega_i=-1$ for all $i\in F_-$, $\omega_i=0$ for all $i\in F_+$, $\omega_i=\mu_i$ for all $i\in I\setminus\{i'\}$, and $\omega_{i'}=\mu_{i'}+s^+(\mu_{F_+}) + s^-(\mu_{F_-})$.
\item $I = \emptyset$ (in particular $|[0]|=1$ and $S(\prec_\mu)=\{\omega(\prec)\}$ where $\omega(\prec)_i=-1$ if $i\prec 0$, and $\omega(\prec)_i=0$ if $i\succ 0$), $s^+(\mu_{F_-})+s^-(\mu_{F_-})=0,1,-1$ and $\mu \sim_{\fsp} \omega(\prec)$, $\mu\sim_{\fsp} \omega(\prec) + \varepsilon_{i'}$, or $\mu\sim_{\fsp} \omega(\prec) - \varepsilon_{i'}$, for some $i'\in \Int^\pm (\mu)$, respectively. Moreover,
 	\[
X_{\fsp}(\mu)^{\fu_{X_\fsp(\mu)}} \cong \C
	\]
as $\fl_{X_\fsp(\mu)}^{ss}=\fsl(F)$-modules.
\end{enumerate}
\end{prop}

\begin{proof}
It follows from Corollary~\ref{cor:X(sp).top}.
\end{proof}

\begin{rem}
Differently from the case of $\fsl(\infty)$ (see Remark~\ref{rem:nice.prop}), it is not true in general that $X_\fsp(\mu)$ is always locally $\p_{X_\fsp(\mu)}$-aligned.
\end{rem}

\begin{lem}
The module $X_{\fsp}(\mu^n)$ is $\p_{X_\fsp(\mu^n)}$-aligned if and only if  $I\neq \emptyset$, or $I=\emptyset$ and $s_n^+(\mu^n_{F_+})+s_n^-(\mu^n_{F_-})\in \{0,\pm 1\}$.
\end{lem}

\begin{proof}
If $i'\in I$, then $\mu^n\sim_{\fsp} \omega$, where $\omega_i=-1$ for every $i\in F_-$, $\omega_i=0$ for every $i\in F_+$, $\omega_i=\mu_i$ for every $i\in I\setminus\{i'\}$, and $\omega_{i'}=\mu_{i'} + s_n^+(\mu^n_{F_+})+s_n^-(\mu^n_{F_-})$. Since $\omega$ is $\fu_{X_\fsp(\mu^n)}$-singular, one direction follows. If $I=\emptyset$, and $s_n^+(\mu^n_{F_+})+s_n^-(\mu^n_{F_-})\in \{0,\pm 1\}$, then the result follows from Corollary~\ref{cor:X_sp.h.w.I=empty}. On the other hand, it is easy to see that if $s_n^+(\mu^n_{F_+})+s_n^-(\mu^n_{F_-})\notin \{0,\pm 1\}$, then either: (1) $s_n^+(\mu^n_{F_+})+s_n^-(\mu^n_{F_-})<-1$, and for any $\lambda\sim_{\fsp}\mu$ we have that either there is $i\in F_-$ such that $\lambda_i<-2$, or there are $i,j\in F_-$ such that $\lambda_i=-2$ and $\lambda_j=-2$. Either way it is clear that ${\bf x}^\lambda$ cannot be $\fu_{X_\fsp(\mu^n)}$-singular; or (2) $s_n^+(\mu^n_{F_+})+s_n^-(\mu^n_{F_-})>1$ and similarly we prove that there is no $\fu_{X_\fsp(\mu^n)}$-singular weight. This along with Corollary~\ref{cor:p-ind.iff.inv.nonzero} proves the result.
\end{proof}

Let $\mu\in \C^{\Z_{>0}}$ and $\mu^n=(\mu_1,\ldots, \mu_n)\in \C^n$. A $\Z_2$-linear partial order $\prec_n$ on $\pm \Z_n\cup\{0\}$ is said to be $\fsp$-\emph{compatible} with $\mu^n$ if $\prec_n$ coincides with $\prec_{\mu^n}$ on $[0]$ and $[0_-]\prec_n [0] \prec_n [0_+]$. A $\Z_2$-linear partial order $\prec$ on $\Z$ is said to be \emph{locally $\fsp$-compatible} with $\mu$ if $\prec_n$ is $\fsp$-compatible with $\mu^n$ for every $n\in \Z_{>0}$.

\begin{lem}\label{lem:(sp).p-aligned.implies.p_mu.aligned}
Consider the above notation. If $X_{\fsp}(\mu^n)$ is $\p(\prec_n)$-aligned, then it is also $\p_{X_\fsp(\mu^n)}$-aligned.
\end{lem}

\begin{proof}
If $I\neq \emptyset$, then $X_\fsp(\mu)$ is always $\p_{X_\fsp(\mu^n)}$-aligned and there is nothing to prove. On the other hand, if $I=\emptyset$ and $X_{\fsp}(\mu^n)$ is $\p(\prec_n)$-aligned, then we claim that $\prec_n$ must be a refinement of $\prec_{\mu^n}$. Indeed, $\prec_{\mu^n}$ induces a decomposition $\pm \Z_n\cup\{0\}=[0_-]\prec_{\mu^n}[0]\prec_{\mu^n} [0_+]$, and it is clear that if $X_{\fsp}(\mu^n)$ is $\p(\prec_n)$-aligned, then we must have $[0_-]\prec_n [0]\prec_n [0_+]$. Thus the claim is proved and the result follows.
\end{proof}

\begin{prop}
Let $\prec$ be a $\Z_2$-linear partial order on $\Z$. Then $X_{\fsp}(\mu)$ is pseudo $\p(\prec)$-aligned if and only if $X_\fsp(\mu)$ is locally $\p_{X_\fsp(\mu)}$-aligned, $\prec$ is locally $\fsp$-compatible with $\mu$, and $X_{\fsp}(\mu)$ is not $\p(\prec)$-aligned.
\end{prop}

\begin{proof}
Notice first that $X_{\fsp}(\mu)\cong \varinjlim_n X_{\fsp}(\mu^n)$. Next, we claim that $X_{\fsp}(\mu^n)$ is $\p(\prec_n)$-aligned if and only if it is $\p_{X_\fsp(\mu^n)}$-aligned and $\prec_n$ is $\fsp$-compatible with $\mu^n$. Indeed, if $I\neq \emptyset$, then $X_\fsp(\mu)$ is $\p_{X_\fsp(\mu^n)}$-aligned for every $n\gg 0$, and we have that $\mu^n\sim_{\fsp} \omega$, where $\omega_i=-1$ for all $i\in \Int^-(\mu^n)$, and $\omega_i=0$ for all $i\in \Int^+(\mu^n)$. If $\prec_n$ is $\fsp$-compatible with $\mu^n$, then 
\begin{enumerate}
\item $i\prec_n 0$ implies $i\in \Int^-(\mu^n)$. Then $\omega_i=-1$ and ${\bf x}^\omega$ is $\g_{2\varepsilon_i}$-singular.
\item $-i\prec_n 0$ implies $i\in \Int^+(\mu^n)$. Then $\omega_i=0$ and ${\bf x}^\omega$ is $\g_{-2\varepsilon_i}$-singular.
\item $i\prec_n j$ implies either $i,j\in \Int^-(\mu^n)$, or $i,j\in \Int^+(\mu^n)$, or $i\in \Int^-(\mu^n)$ and $j\in I\cup \Int^+(\mu^n)$, or $i\in I$ and $j\in \Int^+(\mu^n)$. Then, for all cases, either $\omega_i=-1$ or $\omega_j=0$ and ${\bf x}^\omega$ is $\g_{\varepsilon_i-\varepsilon_j}$-singular.
\item $i\prec_n -j$ implies either $i\in \Int^-(\mu^n)$ and $j\in I\cup \Int^+(\mu^n)\cup \Int^-(\mu^n)$, or $i\in I$ and $j\in  \Int^-(\mu^n)$, or $i\in \Int^+(\mu^n)$ and $j\in \Int^-(\mu^n)$. Then, for all cases, either $\omega_i=-1$ or $\omega_j=-1$ and ${\bf x}^\omega$ is $\g_{\varepsilon_i+\varepsilon_j}$-singular.
\item $-i\prec_n j$ implies either $i\in \Int^+(\mu^n)$ and $j\in I\cup \Int^+(\mu^n)\cup \Int^-(\mu^n)$, or $i\in I$ and $j\in  \Int^+(\mu^n)$, or $i\in \Int^-(\mu^n)$ and $j\in \Int^+(\mu^n)$. Then, for all cases, either $\omega_i=0$ or $\omega_j=0$ and ${\bf x}^\omega$ is $\g_{-\varepsilon_i-\varepsilon_j}$-singular.
\end{enumerate}
Thus ${\bf x}^\omega$ is $\fu(\prec_n)$-singular and $X_{\fsp}(\mu^n)$ is $\p(\prec_n)$-aligned by Corollary~\ref{cor:p-ind.iff.inv.nonzero}. On the other hand, if $\prec_n$ is not $\fsp$-compatible with $\mu^n$, then either there are $i,j\in I$ with $i\prec_n j$, which clearly gives us a contradiction, since for any ${\bf x}^\lambda\in X_\fsl(\mu^n)$, we have that $\lambda_i,\lambda_j\notin \Z$, and hence there is no weight which is $\g_{\varepsilon_i-\varepsilon_j}$-singular; or 
\begin{enumerate}
\item $j\prec_n i$, for $i\in \Int^-(\mu^n)$, $j\in I\cup \Int^+(\mu^n)$;
\item $-j\prec_n i$, for $i\in \Int^-(\mu^n)$, $j\in I\cup \Int^-(\mu^n)$;
\item $j\prec_n -i$, for $i\in \Int^+(\mu^n)$, $j\in I\cup \Int^+(\mu^n)$;
\item $-j\prec_n -i$, for $i\in \Int^+(\mu^n)$, $j\in I\cup \Int^-(\mu^n)$;
\item $j\prec_n i$, for $i\in I$, $j\in \Int^+(\mu^n)$;
\item $-j\prec_n i$, for $i\in I$, $j\in \Int^-(\mu^n)$;
\item $j\prec_n -i$, for $i\in I$, $j\in \Int^+(\mu^n)$;
\item $-j\prec_n -i$, for $i\in I$, $j\in \Int^-(\mu^n)$.
\end{enumerate}
For the cases (a) and (e) there is no weight of $X_{\fsp}(\mu^n)$ which is $\g_{\varepsilon_j-\varepsilon_i}$-singular. For the cases (b) and (f) there is no weight of $X_{\fsp}(\mu^n)$ which is $\g_{-\varepsilon_j-\varepsilon_i}$-singular. For the cases (c) and (g) there is no weight of $X_{\fsp}(\mu^n)$ which is $\g_{\varepsilon_j+\varepsilon_i}$-singular. For the cases (d) and (h) there is no weight of $X_{\fsp}(\mu^n)$ which is $\g_{\varepsilon_i-\varepsilon_j}$-singular (notice that we are using the following fact: if $\nu\in \Supp X_{\fsl}(\eta)$, then $\Int^\pm(\nu) = \Int^\pm (\eta)$). Thus, for the case $I\neq \emptyset$, we have shown that $X_\fsp(\mu^n)$ is $\p(\prec_n)$-aligned if and only if $\prec_n$ is $\fsp$-compatible with $\mu^n$. Suppose now that $I=\emptyset$. Hence, by Lemma~\ref{lem:(sp).p-aligned.implies.p_mu.aligned}, if $X_{\fsp}(\mu^n)$ is not $\p(\prec_{\mu^n})$-aligned then it cannot be $\p(\prec_n)$-aligned. Then we may assume that $X_{\fsp}(\mu^n)$ is $\p(\prec_{\mu^n})$-aligned, and similarly to the previous case, we show that $X_\fsp(\mu^n)$ is $\p(\prec_n)$-aligned if and only if $\prec_n$ is $\fsp$-compatible with $\mu^n$. Therefore $X_\fsp(\mu)$ is locally $\p(\prec)$-aligned if and only if $X_\fsp(\mu)$ is locally $\p_{X_\fsp(\mu)}$-aligned, and $\prec$ is locally $\fsp$-compatible with $\mu$. This concludes the proof.
\end{proof}

\section*{Acknowledgements}

This paper has been written during a post-doctoral period at the Jacobs University, Bremen, under supervision of Ivan Penkov. I am grateful to Ivan Penkov for proposing the problem, for all stimulating discussions, and for valuable suggestions. I also thank Jacobs University for its hospitality.

\section*{Funding}

L. C. was supported by the Capes grant (88881.119190/2016-01) and by the PRPq grant ADRC-05/2016.

\end{document}